\catcode'32=9
\magnification=1200

\font\tenpc=cmcsc10

\font\eightrm=cmr8
\font\eighti=cmmi8
\font\eightsy=cmsy8
\font\eightbf=cmbx8
\font\eighttt=cmtt8
\font\eightit=cmti8
\font\eightsl=cmsl8
\font\sixrm=cmr6
\font\sixi=cmmi6
\font\sixsy=cmsy6
\font\sixbf=cmbx6

\skewchar\eighti='177 \skewchar\sixi='177
\skewchar\eightsy='60 \skewchar\sixsy='60

\font\tengoth=eufm10
\font\tenbboard=msbm10
\font\eightgoth=eufm7 at 8pt
\font\eightbboard=msbm7 at 8pt
\font\sevengoth=eufm7
\font\sevenbboard=msbm7
\font\sixgoth=eufm5 at 6 pt
\font\fivegoth=eufm5

\font\tengoth=eufm10
\font\tenbboard=msbm10
\font\eightgoth=eufm7 at 8pt
\font\eightbboard=msbm7 at 8pt
\font\sevengoth=eufm7
\font\sevenbboard=msbm7
\font\sixgoth=eufm5 at 6 pt
\font\fivegoth=eufm5

\newfam\gothfam
\newfam\bboardfam

\catcode`\@=11

\def\raggedbottom{\topskip 10pt plus 36pt
\r@ggedbottomtrue}
\def\pc#1#2|{{\bigf@ntpc #1\penalty
\@MM\hskip\z@skip\smallf@ntpc #2}}

\def\tenpoint{%
  \textfont0=\tenrm \scriptfont0=\sevenrm \scriptscriptfont0=\fiverm
  \def\rm{\fam\z@\tenrm}%
  \textfont1=\teni \scriptfont1=\seveni \scriptscriptfont1=\fivei
  \def\oldstyle{\fam\@ne\teni}%
  \textfont2=\tensy \scriptfont2=\sevensy \scriptscriptfont2=\fivesy
  \textfont\gothfam=\tengoth \scriptfont\gothfam=\sevengoth
  \scriptscriptfont\gothfam=\fivegoth
  \def\goth{\fam\gothfam\tengoth}%
  \textfont\bboardfam=\tenbboard \scriptfont\bboardfam=\sevenbboard
  \scriptscriptfont\bboardfam=\sevenbboard
  \def\bboard{\fam\bboardfam}%
  \textfont\itfam=\tenit
  \def\it{\fam\itfam\tenit}%
  \textfont\slfam=\tensl
  \def\sl{\fam\slfam\tensl}%
  \textfont\bffam=\tenbf \scriptfont\bffam=\sevenbf
  \scriptscriptfont\bffam=\fivebf
  \def\bf{\fam\bffam\tenbf}%
  \textfont\ttfam=\tentt
  \def\tt{\fam\ttfam\tentt}%
  \abovedisplayskip=12pt plus 3pt minus 9pt
  \abovedisplayshortskip=0pt plus 3pt
  \belowdisplayskip=12pt plus 3pt minus 9pt
  \belowdisplayshortskip=7pt plus 3pt minus 4pt
  \smallskipamount=3pt plus 1pt minus 1pt
  \medskipamount=6pt plus 2pt minus 2pt
  \bigskipamount=12pt plus 4pt minus 4pt
  \normalbaselineskip=12pt
  \setbox\strutbox=\hbox{\vrule height8.5pt depth3.5pt width0pt}%
  \let\bigf@ntpc=\tenrm \let\smallf@ntpc=\sevenrm
  \let\petcap=\tenpc
  \normalbaselines\rm}
\def\eightpoint{%
  \textfont0=\eightrm \scriptfont0=\sixrm \scriptscriptfont0=\fiverm
  \def\rm{\fam\z@\eightrm}%
  \textfont1=\eighti \scriptfont1=\sixi \scriptscriptfont1=\fivei
  \def\oldstyle{\fam\@ne\eighti}%
  \textfont2=\eightsy \scriptfont2=\sixsy \scriptscriptfont2=\fivesy
  \textfont\gothfam=\eightgoth \scriptfont\gothfam=\sixgoth
  \scriptscriptfont\gothfam=\fivegoth
  \def\goth{\fam\gothfam\eightgoth}%
  \textfont\bboardfam=\eightbboard \scriptfont\bboardfam=\sevenbboard
  \scriptscriptfont\bboardfam=\sevenbboard
  \def\bboard{\fam\bboardfam}%
  \textfont\itfam=\eightit
  \def\it{\fam\itfam\eightit}%
  \textfont\slfam=\eightsl
  \def\sl{\fam\slfam\eightsl}%
  \textfont\bffam=\eightbf \scriptfont\bffam=\sixbf
  \scriptscriptfont\bffam=\fivebf
  \def\bf{\fam\bffam\eightbf}%
  \textfont\ttfam=\eighttt
  \def\tt{\fam\ttfam\eighttt}%
  \abovedisplayskip=9pt plus 2pt minus 6pt
  \abovedisplayshortskip=0pt plus 2pt
  \belowdisplayskip=9pt plus 2pt minus 6pt
  \belowdisplayshortskip=5pt plus 2pt minus 3pt
  \smallskipamount=2pt plus 1pt minus 1pt
  \medskipamount=4pt plus 2pt minus 1pt
  \bigskipamount=9pt plus 3pt minus 3pt
  \normalbaselineskip=9pt
  \setbox\strutbox=\hbox{\vrule height7pt depth2pt width0pt}%
  \let\bigf@ntpc=\eightrm \let\smallf@ntpc=\sixrm
  \normalbaselines\rm}

\tenpoint

\frenchspacing


\newif\ifpagetitre
\newtoks\auteurcourant \auteurcourant={\hfil}
\newtoks\titrecourant \titrecourant={\hfil}

\def\appeln@te{}
\def\vfootnote#1{\def\@parameter{#1}\insert\footins\bgroup\eightpoint
  \interlinepenalty\interfootnotelinepenalty
  \splittopskip\ht\strutbox 
  \splitmaxdepth\dp\strutbox \floatingpenalty\@MM
  \leftskip\z@skip \rightskip\z@skip
  \ifx\appeln@te\@parameter\indent \else{\noindent #1\ }\fi
  \footstrut\futurelet\next\fo@t}

\pretolerance=500 \tolerance=1000 \brokenpenalty=5000
\newdimen\hmargehaute \hmargehaute=0cm
\newdimen\lpage \lpage=13.3cm
\newdimen\hpage \hpage=20cm
\newdimen\lmargeext \lmargeext=1cm
\hsize=11.25cm
\vsize=18cm
\parskip 0pt
\parindent=12pt

\def\margehaute{\vbox to \hmargehaute{\vss}}%
\def\margebasse{\vss}

\output{\shipout\vbox to \hpage{\margehaute\nointerlineskip
  \corpsdepage\margebasse}
  \advancepageno \global\pagetitrefalse
  \ifnum\outputpenalty>-20000 \else\dosupereject\fi}

\def\corpsdepage{\hbox to \lpage{\hss\pagetexte\hskip\lmargeext}}
\def\pagetexte{\vbox{\makeheadline\pagebody\makefootline}}
\headline={\ifpagetitre\titleheadline \else
  \ifodd\pageno\rightheadline \else\leftheadline\fi\fi}
\def\leftheadline{\eightpoint\hfil\the\auteurcourant\hfil}
\def\rightheadline{\eightpoint\hfil\the\titrecourant\hfil}
\def\titleheadline{\hfill}
\pagetitretrue

\def\footnoterule{\kern-6\p@
  \hrule width 2truein \kern 5.6\p@} 

\def\pd#1#2 {\pc#1#2| }

\def\pointir{\discretionary{.}{}{.\kern.35em---\kern.7em}\nobreak
\hskip 0em plus .3em minus .4em }

\def\abstract#1{\vbox{\eightpoint \pc ABSTRACT|\pointir #1}}

\def\titre#1|{\message{#1}
              \par\vskip 30pt plus 24pt minus 3pt\penalty -1000
              \vskip 0pt plus -24pt minus 3pt\penalty -1000
              \centerline{\bf #1}
              \vskip 5pt
              \penalty 10000 }

\def\section#1|{\par\vskip .3cm
                {\bf #1}\pointir}

\def\ssection#1|{\par\vskip .2cm
                {\it #1}\pointir}

\long\def\th#1|#2\finth{\par\medskip
              {\petcap #1\pointir}{\it #2}\par\smallskip}

\long\def\tha#1|#2\fintha{\par\medskip
                    {\petcap #1.}\par\nobreak{\it #2}\par\smallskip}

\def\rem#1|{\par\medskip
            {{\it #1}\pointir}}

\def\rema#1|{\par\medskip
             {{\it #1.}\par\nobreak }}

\def\article#1|#2|#3|#4|#5|#6|#7|
    {{\leftskip=7mm\noindent
     \hangindent=2mm\hangafter=1
     \llap{[#1]\hskip.35em}{#2}.\quad
     #3, {\sl #4}, vol.\nobreak\ {\bf #5}, {\oldstyle #6},
     p.\nobreak\ #7.\par}}
\def\livre#1|#2|#3|#4|
    {{\leftskip=7mm\noindent
    \hangindent=2mm\hangafter=1
    \llap{[#1]\hskip.35em}{#2}.\quad
    {\sl #3}.\quad #4.\par}}
\def\divers#1|#2|#3|
    {{\leftskip=7mm\noindent
    \hangindent=2mm\hangafter=1
     \llap{[#1]\hskip.35em}{#2}.\quad
     #3.\par}}
\mathchardef\conj="0365

\def\qed{\quad\raise -2pt\hbox{\vrule\vbox to 10pt{\hrule width 4pt
\vfill\hrule}\vrule}}

\def\decale#1|{\par\noindent\hskip 28pt\llap{#1}\kern 5pt}

\catcode`\@=12
\catcode`\@=12

\def\Grille{\setbox13=\vbox to 5mm{\hrule width 110mm\vfill}
\setbox13=\vbox{\offinterlineskip
   \copy13\copy13\copy13\copy13\copy13\copy13\copy13\copy13
   \copy13\copy13\copy13\copy13\box13\hrule width 110mm}
\setbox14=\hbox to 5mm{\vrule height 65mm\hfill}
\setbox14=\hbox{\copy14\copy14\copy14\copy14\copy14\copy14
   \copy14\copy14\copy14\copy14\copy14\copy14\copy14\copy14
   \copy14\copy14\copy14\copy14\copy14\copy14\copy14\copy14\box14}
\ht14=0pt\dp14=0pt\wd14=0pt
\setbox13=\vbox to 0pt{\vss\box13\offinterlineskip\box14}
\wd13=0pt\box13}


\def\fleche(#1,#2)\dir(#3,#4)\long#5{%
\noalign{\nointerlineskip\leftput(#1,#2){\vector(#3,#4){#5}}\nointerlineskip}}


\def\hfl#1#2#3{\smash{\mathop{\hbox to#3{\rightarrowfill}}\limits
^{\scriptstyle#1}_{\scriptstyle#2}}}

\def\gfl#1#2#3{\smash{\mathop{\hbox to#3{\leftarrowfill}}\limits
^{\scriptstyle#1}_{\scriptstyle#2}}}


 \message{`lline' & `vector' macros from LaTeX}
 \catcode`@=11
\def\{{\relax\ifmmode\lbrace\else$\lbrace$\fi}
\def\}{\relax\ifmmode\rbrace\else$\rbrace$\fi}
\def\newcount{\alloc@0\count\countdef\insc@unt}
\def\newdimen{\alloc@1\dimen\dimendef\insc@unt}
\def\newwrite{\alloc@7\write\chardef\sixt@@n}

\newwrite\@unused
\newcount\@tempcnta
\newcount\@tempcntb
\newdimen\@tempdima
\newdimen\@tempdimb
\newbox\@tempboxa

\def\@spaces{\space\space\space\space}
\def\@whilenoop#1{}
\def\@whiledim#1\do #2{\ifdim #1\relax#2\@iwhiledim{#1\relax#2}\fi}
\def\@iwhiledim#1{\ifdim #1\let\@nextwhile=\@iwhiledim
        \else\let\@nextwhile=\@whilenoop\fi\@nextwhile{#1}}
\def\@badlinearg{\@latexerr{Bad \string\line\space or \string\vector
   \space argument}}
\def\@latexerr#1#2{\begingroup
\edef\@tempc{#2}\expandafter\errhelp\expandafter{\@tempc}%
\def\@eha{Your command was ignored.
^^JType \space I <command> <return> \space to replace it
  with another command,^^Jor \space <return> \space to continue without it.}
\def\@ehb{You've lost some text. \space \@ehc}
\def\@ehc{Try typing \space <return>
  \space to proceed.^^JIf that doesn't work, type \space X <return> \space to
  quit.}
\def\@ehd{You're in trouble here.  \space\@ehc}

\typeout{LaTeX error. \space See LaTeX manual for explanation.^^J
 \space\@spaces\@spaces\@spaces Type \space H <return> \space for
 immediate help.}\errmessage{#1}\endgroup}
\def\typeout#1{{\let\protect\string\immediate\write\@unused{#1}}}

\font\tenln    = line10
\font\tenlnw   = linew10

\newdimen\@wholewidth
\newdimen\@halfwidth
\newdimen\unitlength 

\unitlength =1pt


\def\thinlines{\let\@linefnt\tenln \let\@circlefnt\tencirc
  \@wholewidth\fontdimen8\tenln \@halfwidth .5\@wholewidth}
\def\thicklines{\let\@linefnt\tenlnw \let\@circlefnt\tencircw
  \@wholewidth\fontdimen8\tenlnw \@halfwidth .5\@wholewidth}

\def\linethickness#1{\@wholewidth #1\relax \@halfwidth .5\@wholewidth}

\newif\if@negarg

\def\lline(#1,#2)#3{\@xarg #1\relax \@yarg #2\relax
\@linelen=#3\unitlength
\ifnum\@xarg =0 \@vline
  \else \ifnum\@yarg =0 \@hline \else \@sline\fi
\fi}

\def\@sline{\ifnum\@xarg< 0 \@negargtrue \@xarg -\@xarg \@yyarg -\@yarg
  \else \@negargfalse \@yyarg \@yarg \fi
\ifnum \@yyarg >0 \@tempcnta\@yyarg \else \@tempcnta -\@yyarg \fi
\ifnum\@tempcnta>6 \@badlinearg\@tempcnta0 \fi
\setbox\@linechar\hbox{\@linefnt\@getlinechar(\@xarg,\@yyarg)}%
\ifnum \@yarg >0 \let\@upordown\raise \@clnht\z@
   \else\let\@upordown\lower \@clnht \ht\@linechar\fi
\@clnwd=\wd\@linechar
\if@negarg \hskip -\wd\@linechar \def\@tempa{\hskip -2\wd\@linechar}\else
     \let\@tempa\relax \fi
\@whiledim \@clnwd <\@linelen \do
  {\@upordown\@clnht\copy\@linechar
   \@tempa
   \advance\@clnht \ht\@linechar
   \advance\@clnwd \wd\@linechar}%
\advance\@clnht -\ht\@linechar
\advance\@clnwd -\wd\@linechar
\@tempdima\@linelen\advance\@tempdima -\@clnwd
\@tempdimb\@tempdima\advance\@tempdimb -\wd\@linechar
\if@negarg \hskip -\@tempdimb \else \hskip \@tempdimb \fi
\multiply\@tempdima \@m
\@tempcnta \@tempdima \@tempdima \wd\@linechar \divide\@tempcnta \@tempdima
\@tempdima \ht\@linechar \multiply\@tempdima \@tempcnta
\divide\@tempdima \@m
\advance\@clnht \@tempdima
\ifdim \@linelen <\wd\@linechar
   \hskip \wd\@linechar
  \else\@upordown\@clnht\copy\@linechar\fi}

\def\@hline{\ifnum \@xarg <0 \hskip -\@linelen \fi
\vrule height \@halfwidth depth \@halfwidth width \@linelen
\ifnum \@xarg <0 \hskip -\@linelen \fi}

\def\@getlinechar(#1,#2){\@tempcnta#1\relax\multiply\@tempcnta 8
\advance\@tempcnta -9 \ifnum #2>0 \advance\@tempcnta #2\relax\else
\advance\@tempcnta -#2\relax\advance\@tempcnta 64 \fi
\char\@tempcnta}

\def\vector(#1,#2)#3{\@xarg #1\relax \@yarg #2\relax
\@linelen=#3\unitlength
\ifnum\@xarg =0 \@vvector
  \else \ifnum\@yarg =0 \@hvector \else \@svector\fi
\fi}

\def\@hvector{\@hline\hbox to 0pt{\@linefnt
\ifnum \@xarg <0 \@getlarrow(1,0)\hss\else
    \hss\@getrarrow(1,0)\fi}}

\def\@vvector{\ifnum \@yarg <0 \@downvector \else \@upvector \fi}

\def\@svector{\@sline
\@tempcnta\@yarg \ifnum\@tempcnta <0 \@tempcnta=-\@tempcnta\fi
\ifnum\@tempcnta <5
  \hskip -\wd\@linechar
  \@upordown\@clnht \hbox{\@linefnt  \if@negarg
  \@getlarrow(\@xarg,\@yyarg) \else \@getrarrow(\@xarg,\@yyarg) \fi}%
\else\@badlinearg\fi}

\def\@getlarrow(#1,#2){\ifnum #2 =\z@ \@tempcnta='33\else
\@tempcnta=#1\relax\multiply\@tempcnta \sixt@@n \advance\@tempcnta
-9 \@tempcntb=#2\relax\multiply\@tempcntb \tw@
\ifnum \@tempcntb >0 \advance\@tempcnta \@tempcntb\relax
\else\advance\@tempcnta -\@tempcntb\advance\@tempcnta 64
\fi\fi\char\@tempcnta}

\def\@getrarrow(#1,#2){\@tempcntb=#2\relax
\ifnum\@tempcntb < 0 \@tempcntb=-\@tempcntb\relax\fi
\ifcase \@tempcntb\relax \@tempcnta='55 \or
\ifnum #1<3 \@tempcnta=#1\relax\multiply\@tempcnta
24 \advance\@tempcnta -6 \else \ifnum #1=3 \@tempcnta=49
\else\@tempcnta=58 \fi\fi\or
\ifnum #1<3 \@tempcnta=#1\relax\multiply\@tempcnta
24 \advance\@tempcnta -3 \else \@tempcnta=51\fi\or
\@tempcnta=#1\relax\multiply\@tempcnta
\sixt@@n \advance\@tempcnta -\tw@ \else
\@tempcnta=#1\relax\multiply\@tempcnta
\sixt@@n \advance\@tempcnta 7 \fi\ifnum #2<0 \advance\@tempcnta 64 \fi
\char\@tempcnta}

\def\@vline{\ifnum \@yarg <0 \@downline \else \@upline\fi}

\def\@upline{\hbox to \z@{\hskip -\@halfwidth \vrule
  width \@wholewidth height \@linelen depth \z@\hss}}

\def\@downline{\hbox to \z@{\hskip -\@halfwidth \vrule
  width \@wholewidth height \z@ depth \@linelen \hss}}

\def\@upvector{\@upline\setbox\@tempboxa\hbox{\@linefnt\char'66}\raise
     \@linelen \hbox to\z@{\lower \ht\@tempboxa\box\@tempboxa\hss}}

\def\@downvector{\@downline\lower \@linelen
      \hbox to \z@{\@linefnt\char'77\hss}}

\thinlines

\newcount\@xarg
\newcount\@yarg
\newcount\@yyarg
\newcount\@multicnt
\newdimen\@xdim
\newdimen\@ydim
\newbox\@linechar
\newdimen\@linelen
\newdimen\@clnwd
\newdimen\@clnht
\newdimen\@dashdim
\newbox\@dashbox
\newcount\@dashcnt
 \catcode`@=12


\newbox\tbox
\newbox\tboxa

\def\leftzer#1{\setbox\tbox=\hbox to 0pt{#1\hss}%
     \ht\tbox=0pt \dp\tbox=0pt \box\tbox}

\def\rightzer#1{\setbox\tbox=\hbox to 0pt{\hss #1}%
     \ht\tbox=0pt \dp\tbox=0pt \box\tbox}

\def\centerzer#1{\setbox\tbox=\hbox to 0pt{\hss #1\hss}%
     \ht\tbox=0pt \dp\tbox=0pt \box\tbox}

%
\def\image(#1,#2)#3{\vbox to #1{\offinterlineskip
    \vss #3 \vskip #2}}


\def\leftput(#1,#2)#3{\setbox\tboxa=\hbox{%
    \kern #1\unitlength
    \raise #2\unitlength\hbox{\leftzer{#3}}}%
    \ht\tboxa=0pt \wd\tboxa=0pt \dp\tboxa=0pt\box\tboxa}

\def\rightput(#1,#2)#3{\setbox\tboxa=\hbox{%
    \kern #1\unitlength
    \raise #2\unitlength\hbox{\rightzer{#3}}}%
    \ht\tboxa=0pt \wd\tboxa=0pt \dp\tboxa=0pt\box\tboxa}

\def\centerput(#1,#2)#3{\setbox\tboxa=\hbox{%
    \kern #1\unitlength
    \raise #2\unitlength\hbox{\centerzer{#3}}}%
    \ht\tboxa=0pt \wd\tboxa=0pt \dp\tboxa=0pt\box\tboxa}

\unitlength=1mm

\def\put(#1,#2)#3{\noalign{\nointerlineskip
                               \centerput(#1,#2){$#3$}
                                \nointerlineskip}}
\def\fleche(#1,#2)\dir(#3,#4)\long#5{%
{\leftput(#1,#2){\vector(#3,#4){#5}}}}


\catcode`\@=11
\def\matrice#1{\null \,\vcenter {\normalbaselines \m@th
\ialign {\hfil $##$\hfil &&\  \hfil $##$\hfil\crcr
\mathstrut \crcr \noalign {\kern -\baselineskip } #1\crcr
\mathstrut \crcr \noalign {\kern -\baselineskip }}}\,}

\def\petitematrice#1{\left(\null\vcenter {\normalbaselines \m@th
\ialign {\hfil $##$\hfil 
&&\thinspace  \hfil $##$\hfil\crcr
\mathstrut \crcr \noalign {\kern -\baselineskip } #1\crcr
\mathstrut \crcr \noalign {\kern -\baselineskip }}}\right)}

\catcode`\@=12

\def\inv{\mathop{\rm inv}\nolimits}
\def\imaj{\mathop{\rm imaj}\nolimits}
\def\tot{\mathop{\rm tot}\nolimits}
\def\L{\mathop{\rm L\kern 0pt}}

\def\NIW{\mathop{{\hbox{\eightrm NIW}}}\nolimits}

\def\fix{\mathop{\rm fix}\nolimits}

\def\des{\mathop{\rm des}\nolimits}
\def\dec{\mathop{\rm dec}\nolimits}
\def\maj{\mathop{\rm maj}\nolimits}

\def\pix{\mathop{\rm pix}\nolimits}
\def\wpix{\mathop{\rm wpix}\nolimits}
\def\lec{\mathop{\rm lec}\nolimits}
\def\wlec{\mathop{\rm wlec}\nolimits}
\def\rinv{\mathop{\rm rinv}\nolimits}

\def\Ligne{\mathop{\rm Ligne}\nolimits}
\def\Iligne{\mathop{\rm Iligne}\nolimits}

\def\exc{\mathop{\rm exc}\nolimits}
\def\rmin{\mathop{\rm rmin}\nolimits}
\def\single{\mathop{\rm single}\nolimits}
\def\Single{\mathop{\rm Single}\nolimits}
\def\Lyndon{\mathop{\rm Lyndon}\nolimits}
\def\cyc{\mathop{\rm cyc}\nolimits}
\def\pix{\mathop{\rm pix}\nolimits}
\def\ides{\mathop{\rm ides}\nolimits}
\def\iexc{\mathop{\rm iexc}\nolimits}

\def\lec{\mathop{\rm lec}\nolimits}
\def\dez{\mathop{\rm dez}\nolimits}
\def\maz{\mathop{\rm maz}\nolimits}
\def\maf{\mathop{\rm maf}\nolimits}
\def\LAC{\mathop{{\hbox{\eightrm LAC}}}\nolimits} 
\def\ILAC{\mathop{{\hbox{\eightrm ILAC}}}\nolimits} 
\def\Sym{{\goth S}}

\titrecourant={A QUADRUPLE DISTRIBUTION}
\auteurcourant={DOMINIQUE FOATA AND GUO-NIU HAN}

\rightline{2007/03/12}
\vglue 2cm

\centerline{\bf FIX-MAHONIAN CALCULUS III;}
\smallskip
\centerline{\bf A QUADRUPLE DISTRIBUTION}
\bigskip
\centerline{
\bf Dominique Foata and Guo-Niu Han}

\bigskip\bigskip
\centerline{\bf Abstract}
\medskip
{\narrower\eightpoint

\noindent
A four-variable distribution on
permutations is derived, with two dual combinatorial
interpretations. The first one includes the number of fixed
points ``fix", the second the so-called ``pix" statistic. This
shows that the duality between derangements and
desarrangements can be extended to the case of multivariable
statistics. Several specializations are obtained, including the
joint distribution of (des, exc), where ``des" and ``exc" stand for
the number of descents and excedances, respectively.

}

\bigskip

\centerline{\bf 1. Introduction}

\medskip
Let 
$$\eqalignno{
(a;q)_n&:=\cases{1,&if $n=0$;\cr
(1-a)(1-aq)\cdots (1-aq^{n-1}),&if $n\ge 1$;\cr}\cr
(a;q)_{\infty}&:=\prod_{n\ge 0}(1-aq^n),\cr}
$$
be the traditional notation for the $q$-ascending
factorial. For each $r\ge 0$ 
form the rational fraction
$$
C(r;u,s,q,Y):={(1-qs)\,(u;q)_r\,(usq;q)_r
\over ((u;q)_r-sq(usq;q)_r)(uY;q)_{r+1}}\leqno(1.1)
$$
in four variables $u$, $s$, $q$, $Y$ and expand it as a formal 
power series in~$u$:
$$
C(r;u,s,q,Y)=\sum_{n\ge 0}u^nC_{n}(r;s,q,Y).\leqno(1.2)
$$
It can be verified that each coefficient $C_{n}(r;s,q,Y)$ is actually a
polynomial in three variables with nonnegative integral coefficients.  For
$r,n\ge 0$ consider the set $W_{n}(r)=[0,r]^n$ of all finite words of
length~$n$, whose letters are taken from the alphabet
$[0,r]=\{0,1,\ldots,r\}$. The first purpose of this paper is to show that
$C_{n}(r;s,q,Y)$ is the generating polynomial for $W_{n}(r)$ by two
three-variable statistics $(\dec,\tot,\single)$ and $(\wlec,\tot,\wpix)$,
respectively, defined by means of two classical {\it word factorizations},
the {\it Lyndon factorization} and the $H$-{\it factorization}. See
Theorems 2.1 and 2.3 thereafter and their Corollaries.

\goodbreak
The second purpose of this paper is to consider the formal power series
$$\leqalignno{\qquad
\sum_{r\ge 0}t^r{(1-qs)\,(u;q)_r\,(usq;q)_r
\over ((u;q)_r-sq(usq;q)_r)(uY;q)_{r+1}}
&=\sum_{r\ge 0}t^rC(r;u,s,q,Y)&(1.3)\cr
&=\sum_{r\ge 0}t^r\sum_{n\ge 0}u^nC_{n}(r;s,q,Y),\cr}
$$
expand it as a formal power series in~$u$, but normalized by
denominators of the form $(t;q)_{n+1}$, that is,
$$
(1.4)\ \sum_{r\ge 0}t^r{(1-qs)\,(u;q)_r\,(usq;q)_r
\over ((u;q)_r-sq(usq;q)_r)(uY;q)_{r+1}}\!=\!\sum_{n\ge 0} A_n(s,t,q,Y){u^n\over (t;q)_{n+1}},
$$
and show that each $A_n(s,t,q,Y)$ is actually the {\it generating
polynomial} for the symmetric group~${\goth S}_n$ by two four-variable
statistics $(\exc,\des,\maj,\fix)$ and $(\lec,\ides,\imaj,\pix)$,
respectively. The first (resp. second) statistic involves the number of fixed
points ``fix" (resp. the variable ``pix") and is referred to as the {\it
fix-version} (resp. the {\it pix-version}). Several specializations  of the
polynomials $A_n(s,t,q,Y)$ are then derived with their combinatorial
interpretations. In particular, the {\it joint} distribution of the two
classical Eulerian statistics ``des" and ``exc" is explicitly calculated.

The {\it fix-version} statistic on ${\goth S}_n$, denoted by
$(\exc,\des,\maj,\fix)$, contains the following classical integral-valued
statistics:  the {\it number of excedances} ``exc," the {\it number of
descents} ``des," the {\it major index} ``maj," the {\it number of fixed
points} ``fix," defined for each permutation
$\sigma=\sigma(1)\sigma(2)\cdots
\sigma(n)$ from ${\goth S}_n$ by
$$\eqalignno{\noalign{\vskip-8pt}
\exc\sigma&:=\#\{i:1\le i\le
n-1,\,\sigma(i)>i\};\cr
\des\sigma&:=\#\{i:1\le i\le n-1,\,\sigma(i)>\sigma(i+1)\};\cr
\maj\sigma&:=\sum_i i\quad(1\le i\le
n-1,\,\sigma(i)>\sigma(i+1));\cr
\fix\sigma&:=\#\{i:1\le i\le n,\,\sigma(i)=i\}.\cr}
$$

As was introduced by D\'esarm\'enien {[De84]}, a {\it desarrangement}  is
defined to be a word $w=x_1x_2\cdots x_n$, whose letters are {\it
distinct} positive integers such that the inequalities $x_1>x_2>\cdots
>x_{2j}$ and $x_{2j}<x_{2j+1}$ hold for some $j$ with $1\le j\le n/2$
(by convention: $x_{n+1}=+\infty$). There is no desarrangement of
length~1. Each desarrangement $w=x_1x_2\cdots x_n$ is called a {\it
hook}, if $x_{1}>x_{2}$ and either $n=2$, or $n\ge 3$ and
$x_{2}<x_{3}<\cdots <x_{n}$. As proved by Gessel {[Ge91]}, each
permutation~$\sigma=\sigma(1)\sigma(2)\cdots\sigma(n)$ admits a
unique factorization, called its {\it hook factorization},
$p\tau_1\tau_{2}\cdots \tau_{k}$, where $p$ is an {\it increasing} word
and each factor $\tau_{1}$, $\tau_{2}$, \dots~, $\tau_{k}$ is a hook. To
derive the hook factorization of a permutation, it suffices to start from
the right and at each step determine the right factor which is a hook, or
equivalently, the shortest right factor which is a desarrangement. 

The {\it pix-version} statistic is denoted by
$(\lec,\ides,\imaj,\pix)$. The second and third components are classical: if 
$\sigma^{-1}$ denotes the inverse of the permutation~$\sigma$, they are simply defined by
$$
\eqalignno{\ides\sigma&:=\des\sigma^{-1};\cr
\imaj\sigma&:=\maj\sigma^{-1}.\cr}
$$
The first and fourth components refer to the hook factorization 
$p\tau_1\tau_{2}\cdots \tau_{k}$ of~$\sigma$. For each~$i$ let $\inv\tau_{i}$ denote the {\it number of inversions} of~$\tau_{i}$. Then, we define:
$$
\leqalignno{\lec\sigma&:=\sum_{1\le i\le k}\inv\tau_{i};\cr
\pix\sigma&:=\hbox{length of the factor }p.\cr
}
$$

For instance, the hook factorization of the following permutation of order~14 is indicated by vertical bars. 
$$
\sigma=
1\;3\;4\;14\mid 12\;2\;5\;11\;15\mid 8\;6\;7\mid 13\;9\;10
$$
We have $p=1\;3\;4\;14$, so that $\pix\sigma=4$. Also
$\inv (12\;2\;5\;11\;15)=3$, $\inv(8\;6\;7)=2$, $\inv(13\;9\;10)=2$, so that
$\lec\sigma=7$.
Our main two theorems are the following.

\proclaim Theorem 1.1 {\rm (The fix-version)}. Let 
$A_{n}(s,t,q,Y)$ $(n\ge 0)$ be the sequence of polynomials in four variables, whose factorial generating function is given by $(1.4)$. Then, the generating polynomial for~${\goth S}_{n}$ by the four-variable statistic
$(\exc,\des,\maj,\fix)$ is equal to $A_{n}(s,t,q,Y)$. In other words,
$$\sum_{\sigma\in{\goth S}_n}
s^{\exc\sigma}t^{\des\sigma}q^{\maj\sigma}
Y^{\fix\sigma}=A_n(s,t,q,Y).\leqno(1.5)
$$

\proclaim Theorem 1.2 {\rm (The pix-version)}. Let 
$A_{n}(s,t,q,Y)$ $(n\ge 0)$ be the sequence of polynomials in four variables, whose factorial generating function is given by $(1.4)$. Then, the generating polynomial for ${\goth S}_{n}$ by the four-variable statistic
$(\lec,\ides,\imaj,\pix)$ is equal to $A_{n}(s,t,q,Y)$. In other words,
$$\sum_{\sigma\in{\goth S}_n}
s^{\lec\sigma}t^{\ides\sigma}q^{\imaj\sigma}
Y^{\pix\sigma}=A_n(s,t,q,Y).\leqno(1.6)
$$

The {\it ligne of route}, $\Ligne\sigma$, of a permutation
$\sigma=\sigma(1)\sigma(2)\cdots\sigma(n)$ (also called {\it descent set\/}) is defined to be the {\it set} of all~$i$ such that $1\le i\le n-1$ and $\sigma(i)>\sigma(i+1)$. In particular, $\des\sigma=\#\Ligne\sigma$ and $\maj\sigma$ is the sum of all~$i$ such that $i\in\Ligne\sigma$. Also, let the {\it inverse ligne of route} of~$\sigma$ be defined by $\Iligne\sigma:=\Ligne\sigma^{-1}$, so that $\ides\sigma=\#\Iligne\sigma$ and $\imaj\sigma=\sum_{i}i$ $(i\in\Iligne\sigma)$.
Finally, let $\iexc\sigma:=\exc\sigma^{-1}$.

It follows from Theorem 1.1 and Theorem 1.2 that the two four-variable statistics 
$(\iexc, \ides, \imaj, \fix)$ and $(\lec, \ides, \imaj, \pix)$ are equidistributed on each symmetric group $\Sym_n$. The third goal of this paper is to prove the following stronger result.

\proclaim Theorem 1.3. The two three-variable statistics
$$
(\iexc, \fix, \Iligne)\quad {\sl and}\quad (\lec, \pix, \Iligne)
$$
are equidistributed on each symmetric group~$\Sym_{n}$.

Note that the third component in each of the previous triplets
is a {\it set-valued} statistic. So far, it was known that the
two pairs $(\fix,\Iligne)$ and $(\pix,\Iligne)$ were
equidistributed, a result derived by D\'esarm\'enien and Wachs
{[DeWa88, DeWa93]}, so that Theorem~1.3 may be regarded as an
extension of their result. 
In the following table we reproduce 
the nine derangements (resp. desarrangements)~$\sigma$ from
${\goth S}_4$, which are such that $\fix\sigma=0$ (resp.
$\pix\sigma=0$), together with the values of the pairs
$(\iexc\sigma,\Iligne\sigma)$ (resp.
$(\lec\sigma,\Iligne\sigma)$).

$$
\vbox{\offinterlineskip\halign{\strut
\vrule\ \hfil$#$\hfil\ &
\vrule\ \hfil$#$\hfil\ &
\vrule\ \hfil$#$\hfil\ &
\vrule\ \hfil$#$\hfil\ &
\vrule\ \hfil$#$\hfil\ &
\vrule\ \hfil$#$\hfil\ \vrule\cr
\noalign{\hrule}
\lec&\Iligne&\hbox{\eightpoint Desarrangements }
&\hbox{\eightpoint Derangements}&%
\Iligne&\iexc\cr
\noalign{\hrule}
1&1&2\,1\,3\,4&2\,3\,4\,1&1&1\cr
\noalign{\hrule}
&1,2&3\,2\,4\,1&3\,4\,2\,1&1,2&\cr
\omit\vrule\hfil&\omit\hrulefill&\omit\hrulefill&\omit\hrulefill
&\omit\hrulefill&\omit\hfil\vrule\cr
&1,3&4\,2\,3\,1&2\,4\,1\,3&1,3&\cr
\omit\vrule\hfil&\omit\hrulefill&\omit\hrulefill&\omit\hrulefill
&\omit\hrulefill&\omit\hfil\vrule\cr
\lower 6pt\hbox{2} &\lower 6pt\hbox{2}&3\,1\,2\,4&3\,1\,4\,2&\lower
6pt\hbox{2}&\lower 6pt\hbox{2}\cr
\noalign{\vskip-4pt}
&&3\,1\,4\,2&3\,4\,1\,2&&\cr
\omit\vrule\hfil&\omit\hrulefill&\omit\hrulefill&\omit\hrulefill
&\omit\hrulefill&\omit\hfil\vrule\cr
&1,3&2\,1\,4\,3&2\,1\,4\,3&1,3&\cr
\omit\vrule\hfil&\omit\hrulefill&\omit\hrulefill&\omit\hrulefill
&\omit\hrulefill&\omit\hfil\vrule\cr
&2,3&4\,1\,3\,2&4\,3\,1\,2&2,3&\cr
\omit\vrule\hfil&\omit\hrulefill&\omit\hrulefill&\omit\hrulefill
&\omit\hrulefill&\omit\hfil\vrule\cr
&1,2,3&4\,3\,2\,1&4\,3\,2\,1&1,2,3&\cr
\noalign{\hrule}
3&3&4\,1\,2\,3&4\,1\,2\,3&3&3\cr
\noalign{\hrule}
}}$$

\medskip
In our previous papers {[FoHa06a, FoHa06b]} we have introduced three statistics ``dez," ``maz" and ``maf" on~${\goth S}_{n}$. If $\sigma$ is a permutation, let $i_{1},i_{2},\ldots,i_{h}$ be the increasing sequence of its fixed points. Let $D\sigma$ (resp. $Z\sigma$) be the word derived from~$\sigma=\sigma(1)\sigma(2)\cdots\sigma(n)$ by {\it deleting} all the fixed points (resp. by {\it replacing} all those fixed points by~0). Then those three statistics are simply defined by: $\dez\sigma:=\des Z\sigma$, $\maz\sigma:=\maj Z\sigma$ and
$\maf\sigma:=(i_{1}-1)+(i_{2}-2)+\cdots+(i_{j}-h)+\maj D\sigma$. For instance, with $\sigma=8\,2\,1\,3\,5\,6\,4\,9\,7$ we have $(i_{1},\ldots,i_{h})=(2,5,6)$,
$Z\sigma=8\,0\,1\,3\,0\,0\,4\,9\,7$, $D\sigma=8\,1\,3\,4\,9\,7$ and 
$\dez\sigma=3$, $\maz\sigma=1+4+8=13$, $\maf\sigma=(2-1)+(5-2)+(6-3)+\maj(8\,1\,3\,4\,9\,7)=13$. Theorem 1.4 in [FoHa06a] and Theorem 1.1 above provide 
another combinatorial interpretation for $A_{n}(s,t,q,Y)$, namely
$$
\sum_{\sigma\in {\goth S}_{n}}s^{\exc\sigma}t^{\dez\sigma}
q^{\maz\sigma}Y^{\fix\sigma}=A_{n}(s,t,q,Y).\leqno(1.7)
$$

In the sequel we need the notations for the $q$-multinomial coefficients
$$\displaylines{
{\,n\,\brack
m_1,\ldots,m_k}_q:={(q;q)_n\over (q;q)_{m_1}\cdots
(q;q)_{m_k}}
\quad (m_1+\cdots+m_k= n);\cr
\noalign{\hbox{and the first $q$-exponential}}
e_q(u)=\sum_{n\ge 0}{u^n\over (q;q)_n}={1\over
(u;q)_\infty}.\cr}
$$

\goodbreak
\noindent
Multiply both sides of (1.4) by $1-t$ and
let $t=1$. We obtain the factorial generating function for a sequence of polynomials $(A_n(s,1,q,Y))$ $(n\ge 0)$ in three variables:
$$
\sum_{n\ge 0}
A_n(s,1,q,Y){u^n\over
(q;q)_n} ={(1-sq)e_q(Yu)\over e_q(squ)-sqe_q(u)}.\leqno(1.8)
$$
It follows from Theorem 1.1 that 
$$
\sum_{\sigma\in {\goth S}_n}s^{\exc\sigma}
q^{\maj\sigma}Y^{\fix\sigma}=A_n(s,1,q,Y)\leqno(1.9)$$
holds for every $n\ge 0$, a result stated and proved
by Shareshian and Wachs {[ShWa06]} by means of a symmetric function argument,
so that identity (1.8) with the (1.9) interpretation belongs to those two authors. 
Identity (1.4) can be regarded as a graded form of (1.8). The interest of the graded form also lies in the fact that 
it provides the joint distribution of $(\exc,\des)$,  as shown in (1.15)
below. 

\medskip
Of course, Theorem 1.2 yields a second combinatorial interpretation for the polynomials $A_{n}(s,1,q,Y)$ in the form 
$$
\sum_{\sigma\in {\goth S}_n}s^{\lec\sigma}
q^{\imaj\sigma}Y^{\pix\sigma}=A_n(s,1,q,Y).\leqno(1.10)$$
However we have a third combinatorial interpretation, where the statistic ``imaj" is replaced by the number of inversions ``inv." We state it as our fourth main theorem.

\proclaim Theorem 1.4. Let 
$A_{n}(s,1,q,Y)$ $(n\ge 0)$ be the sequence of polynomials in three variables, whose factorial generating function is given by $(1.8)$. Then, the generating polynomial for ${\goth S}_{n}$ by the three-variable statistic
$(\lec,\inv,\pix)$ is equal to $A_{n}(s,1,q,Y)$. In other words,
$$\sum_{\sigma\in{\goth S}_n}
s^{\lec\sigma}q^{\inv\sigma}
Y^{\pix\sigma}=A_n(s,1,q,Y).\leqno(1.11)
$$

Again, Theorem 1.4 in [FoHa06a] and Theorem 1.1 provide a fourth combinatorial interpretation of $A_{n}(s,1,q,Y)$, namely
$$
\sum_{\sigma\in {\goth S}_n}s^{\exc\sigma}
q^{\maf\sigma}Y^{\fix\sigma}=A_n(s,1,q,Y).\leqno(1.12)$$
Note that the statistic ``maf" was introduced and studied in {[CHZ97]}.

Let $s=1$ in identity (1.4). We get:
$$
 \sum_{n\ge 0}A_n(1,t,q,Y){u^n\over
(t;q)_{n+1}} =\sum_{r\ge 0}t^r\Bigl(
1-u\sum_{i=0}^r q^i\Bigr)^{-1} {(u;q)_{r+1}\over
(uY;q)_{r+1}},\leqno(1.13)
$$
so that Theorem 1.1 implies
$$
\sum_{\sigma\in {\goth S}_n}t^{\des\sigma}
q^{\maj\sigma}Y^{\fix\sigma}=A_n(1,t,q,Y),\leqno(1.14)
$$
an identity derived by Gessel and Reutenauer {[GeRe93]}.

Finally, by letting $q=Y:=1$ we get the generating function for polynomials in two variables $A_{n}(s,t,1,1)$ $(n\ge 0)$ in the form
$$
\sum_{n\ge 0}A_{n}(s,t,1,1){u^n\over (1-t)^{n+1}}
=\sum_{r\ge 0}t^r{1-s\over (1-u)^{r+1}(1-us)^{-r}-s(1-u)}.\leqno(1.15)$$
It then follows from Theorems~1.1 and 1.2 that
$$
\sum_{\sigma\in {\goth S}_n}s^{\exc\sigma}t^{\des\sigma}
=\sum_{\sigma\in {\goth S}_n}s^{\lec\sigma}t^{\ides\sigma}
=A_n(s,t,1,1).\leqno(1.16)
$$
As is well-known (see, {\it e.g.}, {[FoSch70]}) ``exc" and ``des" are equally distributed over~${\goth S}_{n}$, their common generating polynomial being the Eulerian polynomial $A_{n}(t):=A_{n}(t,1,1,1)=A_{n}(1,t,1,1)$, which satifies the identity
$$
{A_{n}(t)\over (1-t)^{n+1}}
=\sum_{r\ge 0}t^r(r+1)^n,\leqno(1.17)
$$
easily deduced from (1.15). 

\goodbreak
The polynomials $A_{n}(s,t,1,1)$ do not have any particular symmetries. 
This is perhaps the reason why their generating function has never been
calculated before, to the best of the authors' knowledge. However, with
$q=1$ and $Y=0$ we obtain
$$
\sum_{n\ge 0}A_{n}(s,t,1,0){u^n\over (1-t)^{n+1}}
=\sum_{r\ge 0}t^r{1-s\over (1-us)^{-r}-s(1-u)^{-r}}.\leqno(1.18)$$
The right-hand side is invariant by the change of variables $u\leftarrow
us$, $s\leftarrow s^{-1}$, so that the polynomials $A_{n}(s,t,1,0)$, which
are the generating polynomials for the set of all {\it derangements} by
the pair $(\exc,\des)$, satisfy $A_{n}(s,t,1,0)=s^n\,A_{n}(s^{-1},t,1,0)$. 
This means that $(\exc, \des)$ and $(\iexc, \des)$ are equidistributed on the set of all derangements. There is a stronger combinatorial result that can be derived as follows. 
Let $\bf c$ be the {\it complement of $(n+1)$} and 
$\bf r$ the {\it reverse image}, which map each permutation
$\sigma=\sigma(1)\ldots \sigma(n)$ onto
${\bf c}\,\sigma:=(n+1-\sigma(1))(n+1-\sigma(2))\ldots (n+1-\sigma(n))$
and ${\bf r}\,\sigma:=\sigma(n)\,\ldots\,\sigma(2)\sigma(1) $, respectively.
Then 
$$
(\exc, \fix, \des, \ides)\,\sigma
=(\iexc, \fix, \des, \ides)\,{\bf c}\,{\bf r}\,\sigma.\leqno(1.19)
$$

\goodbreak
The paper is organized as follows. In order to prove that
$C(r;u,s,q,Y)$ is the generating polynomial for $W_{n}(r)$ by
two multivariable statistics, we show in the next section that
it suffices to construct two explicit bijections $\phi^{\fix}$
and $\phi^{\pix}$. The first bijection, defined in Section~3,
is based on the techniques introduced by Kim and Zeng
[KiZe01]. In
particular, we show that the $V$-cycle decomposition
introduced by those two authors, which is attached to each
permutation, can be extended to the case of words. This is the
content of Theorem~3.4, which may be regarded as our fifth
main result.

The second bijection is constructed in Section~4. In Section~5
we complete the proofs of Theorems~1.1 and~1.2. By combining
the two bijections $\phi^{\fix}$ and $\phi^{\pix}$ we obtain
a transformation on words serving to prove that two bivariable
statistics are equidistributed on the same rearrangement
class.  This is done in Section~6, as well as the proof of
Theorem~1.3. Finally, Theorem~1.4 is proved in Section~7 by
means of a new property of the second fundamental
transformation.

\goodbreak
\bigskip
\centerline{\bf 2. Two multivariable generating functions for words}

\medskip
As $1/(u;q)_r=\sum\limits_{n\ge 0}{r+n-1\brack n}_q\,u^n$
(see, {\it e.g.}, [An76, chap.~3]), we may rewrite the fraction
$C(r;u,s,q,Y)=\displaystyle{(1-qs)\,(u;q)_r\,(usq;q)_r
\over ((u;q)_r-sq(usq;q)_r)(uY;q)_{r+1}}$ as
$$\displaylines{(2.1)\quad
C(r;u,s,q,Y)\hfill\cr
\hfill{}=
\Bigl(1-\sum_{n\ge 2}{r+n-1\brack n}_q u^n
((sq)+(sq)^2+\cdots+(sq)^{n-1})\Bigr)^{-1}\kern-2pt {1\over
(uY;q)_{r+1}}.\cr}
$$
If~$c=c_1c_2\cdots c_n$
is a word, whose letters are nonnegative integers, let 
$\lambda\,c:=n$ be the {\it length} of~$c$ and
$\tot c:=c_1+c_2+\cdots +c_n$ the {\it sum} of its letters.
Furthermore, $\NIW_n$ (resp. $\NIW_n(r)$) designates the set
of all {\it monotonic nonincreasing} words 
$c=c_1c_2\cdots c_n$ of length~$n$, whose letters
are nonnegative integers (resp. nonnegative integers at most
equal to~$r$): $c_{1}\ge c_{2}\ge\cdots\ge c_{1}\ge 0$ (resp.
$r\ge c_{1}\ge c_{2}\ge\cdots\ge c_{1}\ge 0$). Also let $\NIW(r)$ be the union of all
$\NIW_n(r)$ for $n\ge 0$. It is
$q$-routine (see, {\it e.g.}, [An76, chap.~3]) to prove 
$${r+n-1\brack n}_q=\sum_{w\in \hbox{\sixrm NIW}_n(r-1)}
q^{\tot w}.
$$
The sum
$\sum\limits_{n\ge 2}{r+n-1\brack n}_q u^n
((sq)+(sq)^2+\cdots+(sq)^{n-1})$ can then be rewritten as
$\sum\limits_{(w,i)}s^iq^{i+\tot w}u^{\lambda w}$, where the
sum is over all pairs $(w,i)$ such that~$w\in\NIW(r-1)$,
$\lambda w\ge 2$ and~$i$ is an integer satisfying
$1\le i\le \lambda w-1$. Let $D(r)$ (resp. $D_{n}(r)$) denote the
set of all those pairs $(w,i)$ (resp. those pairs such that
$\lambda w=n$). Therefore, equation (2.1) can also be expressed as
$$\displaylines{
C(r;u,s,q,Y)=(1-\kern-5pt\sum_{(w,i)\in D(r)}s^iq^{i+\tot w}u^{\lambda w})^{-1}
\sum_{n\ge 0}u^n\sum_{w\in \hbox{\sixrm NIW}_n(r)}
q^{\tot w}Y^{\lambda w},\cr}$$

$$\displaylines{
\noalign{\hbox{and the coefficient $C_{n}(r;s,q,Y)$ of $u^n$ defined in
(1.3) as}}
(2.2)\ 
C_{n}(r;s,q,Y)=
\sum
s^{i_1+\cdots+i_m}\,q^{i_1+\cdots+i_m+\tot w_0
+\tot w_1+\cdots +\tot w_m}\,Y^{\lambda
w_0},\hfill\cr}
$$
the sum being over all sequences
$(w_{0},(w_1,i_1),\ldots,(w_m,i_m))$ such that $w_0\in \NIW(r)$, each of the pairs
$(w_1,i_1)$, \dots~, $(w_m,i_m)$ belongs to $D(r)$,
 and 
$\lambda w_0+\lambda w_1+\cdots+\lambda w_m=n$. Denote
the set of those sequences by $D_n^*(r)$.

\medskip
The next step is to construct two bijections $\phi^{\fix}$ and $\phi^{\pix}$ of $D_n^*(r)$ onto $W_{n}(r)$ enabling us to calculate certain multivariable statistical distributions {\it on words}. As mentioned in the introduction, $\phi^{\fix}$ relates to the algebra of {\it Lyndon words}, first introduced by Chen, Fox and Lyndon {[Ch58]}, popularized in Combinatorics by 
Sch\"utzenberger {[Sch65]} and now set in common usage in Lothaire
{[Lo83]}. The second bijection $\phi^{\pix}$ relates to the less classical $H$-factorization, the analog for words of the hook factorization introduced by Gessel {[Ge91]}.

Let $l=x_1x_2\cdots x_n$ be a nonempty word, whose letters are
nonnegative integers. Then~$l$ is said to be a {\it Lyndon
word}, if either
$n=1$, or if $n\ge 2$ and, with respect the lexicographic
order, the inequality $x_1x_2\cdots x_n >x_ix_{i+1}\cdots
x_nx_1\cdots x_{i-1}$ holds for every~$i$ such that $2\le
i\le n$. When
$n\ge 2$, we always have $x_1\ge x_i$ for all $i=2,\dots,n$
and $x_i>x_{i+1}$ for at least one integer~$i$ ($1\le i\le
n-1)$, so that it makes sense to define the {\it rightmost
minimal letter} of~$l$, denoted by
$\rmin l$, as
the unique letter $x_{i+1}$ satisfying the inequalities 
$x_i>x_{i+1}$, $x_{i+1}\le x_{i+2}\le \cdots \le x_n$.

Let $w$,~$w'$ be two nonempty primitive words (none of them
can be expressed as $v^b$, where~$v$ is a word and $b$ an
integer greater than or equal to~2). We write $w\preceq w'$ if
and only if $w^{b}\le w'^{b}$, with respect to the lexicographic
order, when~$b$ is large enough. As shown for instance in 
[{Lo83}, Theorem 5.1.5] each nonnempty word~$w$, whose letters
are nonnegative integers, can be written uniquely as a product
$l_1l_2\cdots l_k$, where each~$l_i$ is a Lyndon word and
$l_1\preceq l_2\preceq \cdots
\preceq l_k$. Classically, each Lyndon word is defined to be the
minimum within its class of cyclic rearrangements, so that the
sequence $l_1\preceq l_2\preceq \cdots$ is replaced by
$l_1\ge l_2\ge\cdots$ The modification made here is for
convenience.

For instance, the factorization of the following word as a
nondecreasing product of Lyndon words with respect to
``$\preceq$" [in short, {\it Lyndon word factorization}] is
indicated by vertical bars:
$$
w=\mid 2\mid3\,2\,1\,1\mid 3\mid 5\mid 6\,4\,2\,1\,3\,2\,3\mid
6\,6\,3\,1\,6\,6\,2\mid 6\mid.
$$

Now let $w=x_1x_2\cdots x_n$ be an {\it arbitrary} word. We
say that a positive integer~$i$ is a {\it decrease} of~$w$ if
$1\le i\le n-1$ and $x_i\ge x_{i+1}\ge \cdots \ge
x_j>x_{j+1}$ for some~$j$ such that $i\le j\le n-1$. In
particular, $i$ is a decrease if $x_i>x_{i+1}$. Let $\dec(w)$
denote the {\it number of decreases} of~$w$ (different from the
number of descents). We have $\dec(w)=0$ if all 
letters of~$w$ are equal. Also $\dec(w)\ge 1$ if~$w$ is a Lyndon word
having at least two letters. The number of decreases of the
word~$w$ in the previous example is equal to 11. 

Let $l_1l_2\ldots l_k$ be the Lyndon word factorization of a
word~$w$ and let $(l_{i_1},l_{i_2},\ldots,l_{i_h})$ 
$(1\le i_1<i_2<\cdots <i_h\le k)$ be the
sequence of all the {\it one-letter} factors in its Lyndon word
factorization. Form the nonincreasing word $\Single w$ defined by
$\Single w:=l_{i_h}\cdots l_{i_2}l_{i_1}$ and let 
$\single w=h$ be the number of letters of $\Single w$.
In the previous example we have: $\Single w=6\,5\,3\,2$
and $\single w=4$.

\proclaim Theorem 2.1. The map 
$\phi^{\fix}:(w_{0},(w_{1},i_{1}),\ldots, (w_{m},i_{m}))\mapsto w$
of $D_{n}^*(r)$ onto $W_{n}(r)$, defined in Section~$3$, is a bijection having the properties:
$$\displaylines{\noalign{\vskip-2pt}
(2.3)\qquad\eqalign{&i_1+\cdots+i_m=\dec w;\cr
\noalign{\vskip-2pt}
&i_1+\cdots+i_m+\tot w_0
+\tot w_1+\cdots+\tot w_m
=\tot w;\cr
\noalign{\vskip-2pt}
&\lambda w_0=\single w.\cr}\hfill\cr
}
$$

The next Corollary is then a consequence of (2.2) and the above theorem.

\proclaim Corollary 2.2. The sum $C_{n}(r;s,q,Y)$ defined in $(2.2)$
is also equal to
$$
C_{n}(r;s,q,Y)=\sum_w s^{\dec w}q^{\tot w}Y^{\single w},\leqno(2.4)
$$
where the sum is over all words $w\in W_n(r)$.

\medskip
To define the second bijection
$\phi^{\pix}:D_{n}^*(r)\rightarrow W_n(r)$ another class of words is in use. We call them {\it H-words}. They are defined as follows: let
$h=x_1x_2\cdots x_n$ be a word of length $\lambda h\ge 2$, 
whose letters are nonnegative integers. Say that~$h$ is a {\it H-word}, 
if $x_1<x_2$, and either $n=2$, or $n\ge 3$ and $x_2\geq x_3\geq \cdots \geq x_n$. 

Each nonnempty word~$w$, whose letters
are nonnegative integers, can be written uniquely as a product
$uh_1h_2\cdots h_k$, where $u$ is a monotonic
{\it nonincreasing} word (possibly empty) and each~$h_i$ a $H$-word. This factorization is called the {\it $H$-factorization} of~$w$. Unless~$w$ is monotonic nonincreasing, it ends with a $H$-word, so that its $H$-factorization is  obtained by removing that $H$-word and
determining the next rightmost $H$-word. Note the discrepancy between the hook factorization for {\it permutations} mentioned in the introduction and the present $H$-factorization used for {\it words}.

For instance, the $H$-factorization of the following word 
is indicated by vertical bars:
$$
w=\mid 6\,5\,3\,2\mid1\,3\,2\,1\mid 3\,6\,4\mid 1\,2\mid
2\,3\mid 1\,6\,6\,3\mid2\,6\,6\mid.
$$

\goodbreak
Three statistics are now defined that relate to the $H$-factorization $uh_1h_2\cdots h_k$ of each {\it arbitrary} word~$w$. First, let $\wpix(w)$ be the length $\lambda u$ of $u$. Then,
if $\hbox{\bf r}$ denotes the {\it reverse image}, which maps each word  $x_1x_2\ldots x_n$ onto $x_n\ldots x_2x_1$, define the statistic $\wlec(w)$ by
$$
\wlec(w):=\sum_{i=1}^k \rinv(h_i),
$$
where $\rinv(w)= \inv(\hbox{\bf r}(w))$. 
In the previous example, $\wpix w=\lambda(6532)=4$ and
$\wlec w=\inv(1231)+\inv(463)+\inv(21)+\inv(32)+\inv(3661)+\inv(662)=2+2+1+1+3+2=11$. 

\proclaim Theorem 2.3. The map 
$\phi^{\pix}:(w_{0},(w_{1},i_{1}),\ldots, (w_{m},i_{m}))\mapsto w$
of $D_{n}^*(r)$ onto $W_{n}(r)$, defined in Section~$4$, is a bijection having the properties:
$$\displaylines{
(2.5)\qquad\eqalign{&i_1+\cdots+i_m=\wlec w;\cr
&i_1+\cdots+i_m+\tot w_0
+\tot w_1+\cdots+\tot w_m
=\tot w;\cr
&\lambda w_0=\wpix w.\cr}\hfill\cr
}
$$

\proclaim Corollary 2.4. The sum $C_{n}(r;s,q,Y)$ defined in $(2.2)$
is also equal to
$$
C_{n}(r;s,q,Y)=\sum_w s^{\wlec w}q^{\tot w}Y^{\wpix w},
$$
where the sum is over all words $w\in W_n(r)$.

\bigskip
\centerline{\bf 3. The bijection $\phi^{\fix}$}

\medskip
The construction of the bijection $\phi^{\fix}$ of $D_{n}^*(r)$ onto $W_{n}(r)$ proceeds in four steps and involves three subclasses of Lyndon words: the $V$-words, $U$-words and $L$-words. We can say that $V$- and $U$-words are the word analogs of the $V$- and $U$-cycles introduced by Kim and Zeng {[KiZe01]} for permutations. The present construction is directly inspired by their work.

Each word $w=x_1x_2\cdots x_n$ is said to be a $V$-{\it word}
(resp. a $U$-{\it word\/}), 
if it is of length $n\ge 2$ and 
its letters satisfy  the following inequalities 
$$\leqalignno{
x_1\ge x_2&\ge \cdots \ge
x_i>x_{i+1}\ {\rm and}\ 
x_{i+1}\le x_{i+2}\le \cdots \le x_n<x_i,&(3.1)\cr
\noalign{\hbox{(resp.}}
x_1\ge x_2&\ge \cdots \ge
x_i>x_{i+1}\ {\rm and}\ 
x_{i+1}\le x_{i+2}\le \cdots \le x_n<x_1\ )&(3.2)\cr
}$$
for some $i$ such that $1\le i\le n-1$. Note that if (3.1) or
(3.2) holds, then $\dec(w)=i$. Also $\max w\ \hbox{(the maximum letter of }w)=x_1>x_n$. For example, $v= 5\,5\,4\,\underline 1\,1\,2$ is a $V$-word and $u=8\,7\,\underline 5\,7\,7$ is a $U$-word, but not a $V$-word. Their rightmost minimal letters have been underlined.

\goodbreak
Now $w$ is said to be
a $L$-{\it word}, if it is a Lyndon word of length at least
equal to~2 and whenever $x_1=x_i$ for some~$i$ such that
$2\le i\le n$, then $x_1=x_2=\cdots=x_i$.
For instance, $6\,6\,3\,1\,6\,6\,2$ is a Lyndon word, but not a
$L$-word, but $6\,6\,3\,1\,2\,1$ is a $L$-word.

\medskip
Let $V_n(r)$ (resp. $U_n(r)$, resp. $L_n(r)$, $\Lyndon_n(r)$) be
the set of
$V$-words (resp. $U$-words, resp. $L$-words, resp. Lyndon
words), of length~$n$, whose  letters are at most equal to~$r$.
Also, let
$V(r)$ (resp. $U(r)$, resp. $L(r)$, resp. $\Lyndon(r)$) be the
union of the
$V_n(r)$'s (resp. the $U_n(r)$'s, resp. the $L_n(r)$'s,
resp. the $\Lyndon_n(r)$'s) for $n\ge
2$.  Clearly, $V_n(r)\subset U_n(r)\subset
L_n(r)\subset\Lyndon_n(r)$. Parallel to 
$D_n^*(r)$,  whose definition was given in (2.2), we introduce
three sets $V_n^*(r)$,
$U_n^*(r)$, $L_n^*(r)$ of sequences $(w_0,w_1,\ldots,
w_k)$ of words from
$W(r)$ such that $w_0\in \NIW(r)$, $\lambda w_0+\lambda w_1+\cdots+\lambda w_k=n$ and 

(i) for $V_n^*(r)$ the components~$w_i$ $(1\le i\le k)$ 
belong to $V(r)$;

(ii) for $U_n^*(r)$ the components~$w_i$ $(1\le i\le k)$ 
belong to $U(r)$ and are such that: $\rmin w_1<\max w_2$, $\rmin
w_2<\max w_3$,\ \dots~, $\rmin w_{k-1}<\max w_k$;

(iii) for $L_n^*(r)$ the components~$w_i$ $(1\le i\le k)$ 
belong to $L(r)$ and are such that: $\max w_1\le \max w_2\le \cdots\le 
\max w_k$.

\medskip
The first step consists of mapping the set $D_n(r)$ onto
$V_n(r)$. This is made by means of a very simple bijection,
defined as follows: let $w=x_1x_2\cdots x_n$ be a
nonincreasing word and let $(w,i)$ belong to
$D_n(r)$, so that $n\ge 2$ and $1\le i\le n-1$. Let
$$\displaylines{
y_1:=x_1+1,\ y_2:=x_2+1,\  \ldots ,\ y_i:=x_i+1,\cr
y_{i+1}:=x_n,\ y_{i+2}:=x_{n-1},\ \ldots ,\ y_n:=x_{i+1};\cr
v
:=y_1y_2\ldots y_n.\cr}
$$
The following proposition is evident.

\proclaim Proposition 3.1. The mapping $(w,i)\mapsto v$ is a
bijection of $D_n(r)$ onto $V_n(r)$ satisfying
$\dec(v)=i$ and $\tot v=\tot w+i$.

For instance, the image of $(w=4\,4\,3\,2\,1\,1,\, i=3)$ is the
$V$-word $v= 5\,5\,4\,1\,1\,2$ under the above bijection and
$\dec(v)=3$.

\medskip
Let $(w_0,(w_1,i_1),(w_2,i_2),\ldots, (w_k,i_k))$ belong to $D_n^*(r)$
and, using the bijection of Proposition~3.1, let
$(w_1,i_1)\mapsto v_1$, $(w_2,i_2)\mapsto v_2$,
\dots~, 
$(w_k,i_k)\mapsto v_k$. Then
$$
(w_0,(w_1,i_1),(w_2,i_2),\ldots, (w_k,i_k))
\mapsto(w_0,v_1,v_2,\ldots,v_k)\leqno(3.3)$$ is a bijection of
$D_n^*(r)$ onto $V_n^*(r)$ having the property that
$$\eqalign{
&i_1+i_2+\cdots+i_k=\dec(v_1)+\dec(v_2)+\cdots+\dec(v_k);\cr
&i_1+i_2+\cdots+i_k+\tot w_0+\tot w_1+\tot w_2+\cdots+\tot w_k\cr
&\kern100pt{}
=\tot w_0+\tot v_1+\tot v_2+\cdots+\tot v_k.\cr}\leqno(3.4)
$$

\medskip
For instance, the sequence
$$\displaylines{
\bigl(6\,5\,3\,2,\ (2\,1\,1\,1,\ 2), (5\,3\,3,\ 2), (1\,1,\ 1),
(2\,2,\ 1), (5\,5\,2\,1,\ 3), (5\,5\,2,\ 2)\bigr)\cr
\noalign{\hbox{from $D_{22}^*(6)$
is mapped under (3.3)
onto the sequence}}
(6\,5\,3\,2,\ 3\,2\,1\,1,\ 6\,4\,3,\ 2\,1,\ 3\,2,\ 6\,6\,3\,1,\ 6\,6\,2)\in V_{22}^*(6).\cr}
$$
Also $\dec(3\,2\,1\,1)+\dec( 6\,4\,3)+\dec( 2\,1)+\dec(3\,2)+
\dec(6\,6\,3\,1)+\dec( 6\,6\,2)=2+2+1+1+3+2=11$. 

\medskip
The second step is to map $V_n^*(r)$ onto
$U_n^*(r)$. Let $u=y_1y_2\cdots y_k\in
U(r)$ and 
$v=z_1z_2\cdots z_l\in V(r)$. Suppose that
$\rmin u$ is the $(i+1)$-st leftmost letter
of~$u$ and $\rmin v$ is the $(j+1)$-st letter of~$v$.
Also assume that $\rmin u\ge \max v$. Then, the word
$[u,v]:=y_1\cdots y_i z_1\cdots z_jz_{j+1}\cdots z_l
y_{i+1}\cdots y_k$ belongs to $U(r)$. Furthermore, $\rmin
[u,v]$ is the $(i+j+1)$-st leftmost letter of $[u,v]$ and its
value is $z_{j+1}$. We also have the inequalities:
$y_i>y_{i+1}$ (by definition of $\rmin u$),
$y_{i+1}\ge z_1$ (since $\min u\ge \max v$)
and $z_j>z_l$ (since $v$ is a $V$-word). 
These properties allow us 
to get back the pair $(u,v)$ from $[u,v]$ by {\it successively}
determining the {\it critical} letters
$z_{j+1}$, $z_j$, $z_l$, $y_{i+1}$, $z_1$.
The mapping $(u,v)\mapsto [u,v]$ is perfectly reversible.

For example, with $u=8\,7\,\underline5\,7\,7$ and 
$v={\bf 5\,\underline2\,2}$ we have
$[u,v]
=8\,7\,{\bf 5\,\underline 2\,2}\,5 \,7\,7$. 
\smallskip
Now let
$(w_0,v_1,v_2,\ldots, v_k)\in V_n^*(r)$. If $k=1$, let
$(w_0,u_1):=(w_0,v_{1})\in U_n^*(r)$.
If $k\ge 2$, let $(1,2,\ldots,a)$ be the longest sequence of
integers such that
$\rmin v_1\ge \max v_2>\rmin v_2\ge \max v_2>\cdots
\ge\max v_a>\rmin v_a$ and, either $a=k$, or $a\le k-1$ and
$\rmin v_a<\max v_{a+1}$. Let
$$
u_1:=
\cases{v_1,&if $a=1$;\cr
[\cdots[\,[v_1,v_2],v_3],\,\cdots\,, v_a],&if $a\ge 2$.\cr}
\leqno(3.5)$$
We have $u_1\in U(r)$ and $(w_0,u_1)\in U_n^*(r)$ if $a=k$.
Otherwise,
$\rmin u_1<\max v_{a+1}$. We can then
apply the procedure described in (3.5) to the sequence
$(v_{a+1},v_{a+2},\ldots,v_k)$. When reaching~$v_k$ we
obtain a sequence $(w_0,u_1,\ldots,u_h)\in U_n^*(r)$. The
whole procedure is perfectly reversible. We have then the
following proposition.

\proclaim Proposition 3.2. The mapping
$$(w_0,v_1,v_2,\ldots, v_k)\mapsto 
(w_0,u_1,\!u_2,\ldots, u_h)
\leqno(3.6)$$ 
described in $(3.5)$ is a
bijection of
$V^*_n(r)$ onto
$U_n^*(r)$ having the following properties:\hfil\break
(i) $u_1u_2\cdots u_h$ is a rearrangement of
$v_1v_2\cdots v_k$, so that 
$\tot u_1+\tot u_2+\cdots +\tot u_h
=\tot v_1+\tot v_2+\cdots +\tot v_k$;\hfil\break
(ii) $\dec( u_1)+\dec( u_2)+\cdots +\dec( u_h)
=\dec( v_1)+\dec( v_2)+\cdots +\dec( v_k)$.

For example, the above sequence
$$\displaylines{
(6\,5\,3\,2,\ 3\,2\,1\,1,\ 6\,4\,3,\ 2\,1,\ 3\,2,\ 6\,6\,3\,1,\ 6\,6\,2)\in V_{22}^*(6)\cr
\noalign{\hbox{is mapped onto the sequence}}
(6\,5\,3\,2,\ 3\,2\,1\,1,\ 6\,4\,2\,1\,3,\  3\,2,\ 6\,6\,3\,1,\ 6\,6\,2)\in U_{22}^*(6).\cr}
$$
where $[6\,4\,3,\ 2\,1]=6\,4\,2\,1\,3$.
Also $\dec(3\,2\,1\,1)+\dec( 6\,4\,2\,1\,3)+\dec(  3\,2)
+\dec(6\,6\,3\,1)+\dec( 6\,6\,2 )=11$.

\medskip
The third step is to map $U_n^*(r)$ onto $L_n^*(r)$.
Let $l=x_1x_2\cdots x_j\in L(r)$ and
$u=y_1y_2\cdots y_{j'}\in U(r)$. 
Suppose that $\rmin l$ is the
$(i+1)$-st leftmost letter of~$l$ and $\rmin u$ is the
$(i'+1)$-st leftmost letter of~$u$. Also assume that $\rmin l<\max u$ and
$\max l>\max u$. 
If $x_j<y_1$, let
$\langle l,u\rangle:=lu$. If $x_{j}\ge y_{1}$,
there is a unique integer~$a\ge
i+1$ such that
$x_a<y_1=\max u\le x_{a+1}$. Then, let
$$
\langle l,u\rangle:=x_1\cdots x_ix_{i+1}\cdots
x_a y_1\cdots y_{i'}y_{i'+1}\cdots y_{j'}  x_{a+1}\cdots x_j.$$
The word $\langle l,u\rangle$ belongs to $L(r)$ and 
$\rmin \langle l,u\rangle$ is the $(a+i'+1)$-st leftmost letter
of $\langle l,u\rangle$, its value being $y_{i'+1}$. Now $y_1$
is the rightmost letter of $\langle l,u\rangle$ to the left
of~$\rmin\langle l,u\rangle=y_{i'+1}$ such that the letter
preceding it, that is~$x_a$, satisfies
$x_a<y_1$ and the letter following it, that is $y_2$,  is
such that $y_1\ge y_2$. On the other hand, $y_{j'}$ is the
unique letter in the nondecreasing factor $y_{i'+1}\cdots
y_{j'}x_{a+1}\cdots x_j$ that satisfies
$y_{j'}<y_1=\max u\le x_{a+1}$. Hence the mapping
$(l,u)\mapsto\langle l,u\rangle$ is completely reversible.

\medskip
For example, with $l=
8\,2\,5\,\underline
3\,3\,5\,7\,7$ and $u=\bf7\,6\,\underline1\,2$ we have
$\langle l,u\rangle=
8\,2\,5\,
3\,3\,5\,{\bf7\,6\,\underline1\,2}\,7\,7$ (the leftmost minimal
letters have been underlined). The letter $\bf 7$ is the
rightmost letter in $\langle l,u\rangle$ greater than its
predecessor~5 and greater than or equal to its successor~{\bf 6}.
Also {\bf2} is the unique letter in
the factor $1\,{\bf2}\,7\,7$ that satisfies ${\bf 2}<\max
u={\bf7}\le 7$. 

\medskip
Let $(w_0,u_1,u_2,\ldots, u_h)\in U_n^*(r)$. If $h=1$, let 
$(w_0,l_1):=(w_0,u_1)\in L_n^*(r)$. If $h\ge 2$, let
$(1,2,\ldots,a)$ be the longest sequence of integers such that
$\max u_1>\max u_j$ for all $j=2,\ldots, a$ and, and either
$a=h$, or $a\le h-1$ and $\max u_1\le \max u_{a+1}$. Let
$$
l_1:=\cases {u_1,&if $a=1$;\cr
\langle\cdots\langle\,\langle
u_1,u_2\rangle,u_3\rangle,\,\cdots\,,u_a\rangle,&if $a\ge 2$.\cr}
\leqno(3.7)
$$
We have $l_1\in L(r)$ and $(w_0,l_1)\in L_n^*(r)$ if $a=h$.
Otherwise, $\max l_{1}\le \max u_{a+1}$. We then
apply the procedure described in (3.7) to the sequence
$(u_{a+1},u_{a+2},\ldots,u_h)$. When reaching~$u_h$ we
obtain a sequence $(w_0,l_1,\ldots,l_m)\in L_n^*(r)$. The
whole procedure is perfectly reversible. We have then the
following proposition.

\proclaim Proposition 3.3. The mapping
$$
(w_0,u_1,u_2,\ldots,u_h)\mapsto (w_0,l_1,l_2,\ldots,l_m)\leqno(3.8)
$$
described in $(3.7)$ is a bijection of $U_n^*(r)$ onto 
$L_n^*(r)$ having the following properties:\hfil\break
(i) $l_1l_2\cdots l_m$ is a rearrangement of
$u_1u_2\cdots u_h$, so that 
$\tot l_1+\tot l_2+\cdots +\tot l_m=\tot u_1+\tot
u_2+\cdots +\tot u_h $;\hfil\break
(ii) $\dec( l_1)+\dec( l_2)+\cdots +\dec( l_m)=\dec(
u_1)+\dec( u_2)+\cdots +\dec( u_h)$.

For example the above sequence
$$
\displaylines{
(6\,5\,3\,2,\ 3\,2\,1\,1,\ 6\,4\,2\,1\,3,\  3\,2,\ 6\,6\,3\,1,\ 6\,6\,2)
)\in U_{22}^*(6).\cr
\noalign{\hbox{is mapped onto the sequence}}
(6\,5\,3\,2,\ 3\,2\,1\,1,\ 6\,4\,2\,1\,3\,  2\,3,\ 6\,6\,3\,1,\ 6\,6\,2)\in L_{22}^*(6),\cr
}
$$
where $\langle 6\,4\,2\,1\,3,\  3\,2\rangle=
6\,4\,2\,1\,3\,  2\,3$. Also
$\dec(3\,2\,1\,1)+\dec( 6\,4\,2\,1\,3\,  2\,3)+\dec( 6\,6\,3\,1)
+\dec(6\,6\,2)=11$.

\goodbreak
\medskip
The fourth step is to map $L_n^*(r)$ onto $W_n(r)$. Let
$(w_0,l_1,l_2,\ldots,l_m)\in L_n^*(r)$. If $w_0$ is nonempty,
of length~$b$, denote by $f_1$, $f_2$, \dots~, $f_b$ its~$b$
letters from left to right, so that 
$r\ge f_1\ge f_2\ge \cdots\ge f_b\ge 0$. If $m=1$, let
$\sigma_1:=l_1$. If $m\ge 2$, let $a$ be the
greatest integer such that 
$l_1\succ l_2$, $l_1l_2\succ l_3$, \dots~, $l_1\cdots
l_{a-1}\succ l_a$. If $a\le h-1$, let~$a'>a$ be the greatest
integer such that $l_{a+1}\succ l_{a+2}$, $l_{a+1}l_{a+2}\succ
l_{a+3}$, \dots~, $l_{a+1}\cdots l_{a'-1}\succ l_{a'}$, etc.
Form
$\sigma_1:=l_1l_2\cdots l_a$, $\sigma_2:=l_{a+1}\cdots
l_{a'}$, etc. The sequence $(\sigma_1,\sigma_2,\ldots\,)$ is a
{\it nonincreasing} sequence of Lyndon words. Let
$(\tau_1,\tau_2,\ldots,\tau_p)$ be the {\it nonincreasing}
rearrangement of the sequence
$(\sigma_1,\sigma_2,\ldots,f_1,f_2,\ldots,f_b)$ if~$w_0$ is
nonempty, and of $(\sigma_1,\sigma_2,\ldots\,)$ otherwise.
Then, $(\tau_1,\tau_2,\ldots,\tau_p)$ is the {\it Lyndon word
factorization} of a unique word $w\in W_n(r)$. The mapping 
$$(w_0,l_1,l_2,\ldots,l_m)\mapsto w\leqno(3.9)$$ is perfectly reversible. Also
the verification of
$\dec(w)=\dec(l_1)+\dec(l_2)+\cdots+\dec(l_m)$ is immediate.

\medskip
For example, the above sequence
$$
\displaylines{
(6\,5\,3\,2,\ 3\,2\,1\,1,\ 6\,4\,2\,1\,3\,  2\,3,\ 6\,6\,3\,1,\ 6\,6\,2)\in L_{22}^*(6)\cr
\noalign{\hbox{is mapped onto the Lyndon word factorization}}
\rlap{(3.10)}\hfill
w=\mid 2\mid3\,2\,1\,1\mid 3\mid 5\mid 6\,4\,2\,1\,3\,2\,3\mid
6\,6\,3\,1\,6\,6\,2\mid 6\mid.\hfill\cr}
$$
Also $\dec( w)=11$.

\medskip
The map $\phi^{\fix}$ is then defined as being the composition product of
$$
\eqalignno{
(w_0,(w_1,i_1),(w_2,i_2),\ldots, (w_k,i_k))
&\mapsto(w_0,v_1,v_2,\ldots,v_k)\qquad&\hbox{in }(3.3)\cr
(w_0,v_1,v_2,\ldots, v_k)&\mapsto 
(w_0,u_1,\!u_2,\ldots, u_h)
&\hbox{in }(3.6)\cr
(w_0,u_1,u_2,\ldots,u_h)&\mapsto (w_0,l_1,l_2,\ldots,l_m)
&\hbox{in }(3.8)\cr
(w_0,l_1,l_2,\ldots,l_m)&\mapsto w.&\hbox{in }(3.9)\cr
}
$$
Therefore, $\phi^{\fix}$ is a bijection of $D_{n}^*(r)$ onto $W_{n}(r)$ having the properties stated in Theorem~2.1. The latter theorem is then proved.

\medskip
From the property of the bijection $w\mapsto
(w_0,v_1,v_2,\ldots,v_k)$ of $W_{n}(r)$ onto $V_{n}^*(r)$ we
deduce the following theorem.

\proclaim Theorem 3.4 {\rm($V$-word decomposition)}. To each
word
$w=x_{1}x_{2}\cdots x_{n}$ whose letters are nonnegative
integers there corresponds a unique sequence
$(w_0,v_1,v_2,\ldots,v_k)$, where $w_{0}$ is a nondecreasing word and $v_{1}$, $v_{2}$, \dots~, $v_{k}$ are $V$-words with the further property that $w_{0}v_{1}v_{2}\cdots v_{k}$ is a rearrangement of~$w$ and $\dec w=\dec v_{1}+\dec v_{2}+\cdots+\dec v_{k}$.

This theorem may be considered as a {\it word analog}
of Theorem~2.4 in Kim-Zeng's paper {[KiZe01]}. In the
following table the $V$-word decomposition attached to each
rearrangement~$w$ of $1\,2\,2\,3$ has been
calculated, together with the value of $\dec w$. The symbol
``$e$" stands for the empty sequence.

$$
\vbox{\offinterlineskip\halign{\strut\vrule\quad  \hfil$#$\hfil
&\quad \vrule \quad\hfil$#$\hfil\quad&\vrule\quad \hfil$#$\hfil\quad\vrule\cr
\noalign{\hrule}
w&\dec w&(w_{0};v_{1},\ldots,v_{k})\cr
\noalign{\hrule}
1\,2\,2\,3&0&(1\,2\,2\,3;\ e)\cr
\noalign{\hrule}
1\,2\,3\,2&1&(1\,2;\ 3\,2)\cr
\noalign{\hrule}
1\,3\,2\,2&1&(1;\ 3\,2\,2)\cr
\noalign{\hrule}
2\,1\,2\,3&1&(2\,3;\ 2\,1)\cr
\noalign{\hrule}
2\,1\,3\,2&2&(e;\ 2\,1,\ 3\,2)\cr
\noalign{\hrule}
2\,2\,1\,3&2&(3;\ 2\,2\,1)\cr
\noalign{\hrule}
2\,2\,3\,1&1&(2\,2;\ 3\,1)\cr
\noalign{\hrule}
2\,3\,1\,2&1&(2;\ 3\,1\,2)\cr
\noalign{\hrule}
2\,3\,2\,1&2&(2;\ 3\,2\,1)\cr
\noalign{\hrule}
3\,1\,2\,2&1&(e;\ 3\,1\,2\,2)\cr
\noalign{\hrule}
3\,2\,1\,2&2&(e;\ 3\,2,\ 2\,1)\cr
\noalign{\hrule}
3\,2\,2\,1&3&(e;\ 3\,2\,2\,1)\cr
\noalign{\hrule}
}}
$$

\bigskip
\centerline{\bf 4. The bijection $\phi^{\pix}$}

\medskip
In the introduction a hook was defined to be a word $x_1x_2\cdots x_n$
with distinct letters such that $x_1>x_2$ and either $n=2$, or $n\ge3$ and
$x_2<x_3<\cdots <x_n$. In the next definition the letters can be repeated and the inequalities are reversed. Let $h=x_1x_2\cdots x_n$ be a word of length   $\lambda h\ge 2$, whose letters are nonnegative integers. Say that~$h$ is a 
{\it H-word}, if $x_1<x_2$ and either $n=2$, or $n\ge 3$ and $x_2\geq x_3\geq \cdots \geq x_n$.

Each nonnempty word~$w$, whose letters
are nonnegative integers, can be written uniquely as a product
$uh_1h_2\cdots h_k$, where $u$ is a monotonic
{\it nonincreasing} word (possibly empty) and each~$h_i$ is a $H$-word. This factorization is called the {\it $H$-factorization} of~$w$. Unless~$w$ is monotonic nonincreasing, it ends with a $H$-word, so that its $H$-factorization is easily obtained by removing that $H$-word and
determining the next rightmost $H$-word.
For instance, the $H$-factorization of the following word 
is indicated by vertical bars:
$$
w=\mid 6\,5\,3\,2\mid1\,3\,2\,1\mid 3\,6\,4\mid 1\,2\mid
2\,3\mid 1\,6\,6\,3\mid2\,6\,6\mid.
$$

Three statistics are now defined that relate to the $H$-factorization $uh_1h_2\cdots h_k$ of each {\it arbitrary} word~$w$. First, let $\wpix(w)$ be the length $\lambda u$ of~$u$. Then,
if $\hbox{\bf r}$ denotes the {\it reverse image}, which maps each word  $x_1x_2\ldots x_n$ onto $x_n\ldots x_2x_1$, let $\rinv:=\inv\circ\,{\bf r}$ and define the statistic $\wlec(w)$ by
$$
\wlec(w):=\sum_{i=1}^k \rinv(h_i).
$$
In the previous example, $\wpix w=\lambda(6532)=4$ and
$\wlec w=\inv(1231)+\inv(463)+\inv(21)+\inv(32)+\inv(3661)+\inv(662)=2+2+1+1+3+2=11$.

\medskip
The bijection $\phi^{\pix}:D_{n}^*(r)\rightarrow W_{n}(r)$ 
whose properties were stated in Theorem~3.2 is easy to construct.
Let $H_n(r)$  be the set of all
$H$-words of length~$n$, whose  letters are at most equal 
to~$r$ and $H(r)$  be the union of all $H_n(r)$'s  for $n\ge 2$. 
We first map $D_n(r)$ onto $H_n(r)$ as follows.
Let $w=x_1x_2\cdots x_n$ be a
nonincreasing word and let $(w,i)$ belong to
$D_n(r)$, so that $n\ge 2$ and $1\le i\le n-1$. Define:
$$
h :=x_{i+1} (x_1+1)(x_2+1)\ldots (x_i+1) x_{i+2}x_{i+3}\ldots x_n.
$$
The following proposition is evident.

\proclaim Proposition 4.1. The mapping $(w,i)\mapsto h$ is a
bijection of $D_n(r)$ onto $H_n(r)$ satisfying
$\rinv(h)=i$ and $\tot h=\tot w+i$.

For instance, the image of $(w=4\,4\,3\,2\,2\,1,\, i=3)$ is the
$H$-word $h= 2\,5\,5\,4\,2\,1$ under the above bijection and
$\rinv(h)=3$.

\medskip
Let $(w_0, (w_1,i_1),(w_2,i_2),\ldots, (w_k,i_k))$ belong to $D_n^*(r)$
and, using the bijection of Proposition~4.1, let
$(w_1,i_1)\mapsto h_1$, $(w_2,i_2)\mapsto h_2$,
\dots~, 
$(w_k,i_k)\mapsto h_k$. Then $w_0h_1h_{2}\cdots h_k$ is the $H$-factorization of a word $w\in W_{n}(r)$. Accordingly,
$$
\phi^{\pix}:(w_0, (w_1,i_1),(w_2,i_2),\ldots, (w_k,i_k))
\mapsto w:=w_0h_1h_2\cdots h_k$$ 
is a bijection of
$D_n^*(r)$ onto $W_{n}(r)$ having the properties listed in (2.5).
This completes the proof of Theorem~2.3.

\medskip
For instance, the sequence
$$\displaylines{
\bigl(6\,5\,3\,2,\ (2\,1\,1\,1,\ 2), (5\,3\,3,\ 2), (1\,1,\ 1),
(2\,2,\ 1), (5\,5\,2\,1,\ 3), (5\,5\,2,\ 2)\bigr)\cr
\noalign{\hbox{from $D_{22}^*(6)$
is mapped under $\phi^{\pix}$
onto the word}}
6\,5\,3\,2\mid 1\,3\,2\,1\mid 3\,6\,4\mid 1\,2\mid 2\,3\mid 1\,6\,6\,3\mid 2\,6\,6\in W_{22}(6).\cr}
$$
Also $(\wlec,\tot,\wpix)\,w=(11,74,4)$.

\bigskip
\goodbreak
\centerline{\bf 5. From words to permutations}

\medskip
We are now in a position to prove Theorems~1.1 and 1.2.
Suppose that identity (1.5) holds.
As $${1\over(t;q)_{n+1}}=\sum_{j\ge 0}t^j\
\sum\limits_{w\in \hbox{\sevenrm NIW}_n(j)}q^{\tot w},$$
the right-hand side of (1.4) can then be written as
\smash{$\sum\limits_{r\ge 0}t^r\sum\limits_{n\ge 0}B^{\fix}_n(r;s,q,Y)$}, where
$$
B^{\fix}_{n}(r;s,q,Y):=
\sum_{(\sigma,c)} s^{\exc\sigma}\,
q^{\maj\sigma+\tot c}\,Y^{\fix\sigma},\leqno(5.1)
$$
the sum being over all pairs $(\sigma,c)$ such that
$\sigma\in{\goth S}_n$, $\des\sigma\le r$
and $c\in \NIW_n(r-\des\sigma)$. Denote the set of all those pairs by
${\goth S}_{n}(r,\des)$.

In the same manner, let ${\goth S}_{n}(r,\ides)$ denote the set of all pairs $(\sigma,c)$ such that
$\sigma\in{\goth S}_n$, $\ides\sigma\le r$
and $c\in \NIW_n(r-\ides\sigma)$ and let
$$
B^{\pix}_{n}(r;s,q,Y):=
\sum_{(\sigma,c)} s^{\lec\sigma}\,
q^{\imaj\sigma+\tot c}\,Y^{\pix\sigma},\leqno(5.2)
$$
where the sum is over all $(\sigma,c)\in {\goth S}_{n}(r,\ides)$.
If (1.6) holds, the right-hand side of (1.4) is equal to
$\sum\limits_{r\ge 0}t^r\sum\limits_{n\ge 0}B^{\pix}_{n}(r;s,q,Y)$.

Accordingly, for proving identity (1.5) (resp. (1.6)) it suffices to show that
$C_{n}(r;s,q,Y)=B^{\fix}_{n}(r;s,q,Y)$ (resp. $C_{n}(r;s,q,Y)=B^{\pix}_{n}(r;s,q,Y)$) holds for all pairs $(r,n)$. Referring to Corollaries~2.2 and~2.4 it suffices to construct a bijection $$\psi^{\fix}:w\mapsto (\sigma,c)$$ of $W_{n}(r)$ onto
${\goth S}_n(r,\des)$ having the following properties
$$
\eqalign{
\dec w&=\exc \sigma;\cr
\tot w&=\maj\sigma+\tot c;\cr
\single w&=\fix\sigma;\cr}\leqno(5.3)
$$
and a bijection $$\psi^{\pix}:w\mapsto (\sigma,c)$$ of $W_{n}(r)$ onto
${\goth S}_n(r,\ides)$ having the following properties
$$
\eqalign{
\wlec w&=\lec \sigma;\cr
\tot w&=\imaj\sigma+\tot c;\cr
\wpix w&=\pix\sigma.\cr}\leqno(5.4)
$$

The construction of $\psi^{\fix}$ is achieved by adapting a classical
bijection used by Gessel-Reutenauer {[GeRe93]} and
D\'esarm\'enien-Wachs {[DeWa88, DeWa93]}. Start with the
Lyndon word factorization $(\tau_1,\tau_2,\ldots,\tau_p)$ of a
word~$w\in W_n(r)$. If~$x$ is a letter of the
factor $\tau_i=y_1\cdots y_{j-1}xy_{j+1}\cdots y_h$, form
the cyclic rearrangement $\cyc(x):=xy_{j+1}\cdots y_h
y_1\cdots y_{j-1}$. If $x$, $y$ are two letters of~$w$, we say
that~$x$ {\it precedes}~$y$, if $\cyc x\succ \cyc y$, or if 
$\cyc x=\cyc y$ and the letter~$x$ is to the right of the
letter~$y$ in the word~$w$. Accordingly, to each letter~$x$
of~$w$ there corresponds a unique integer~$p(x)$, which is
the number of letters {\it preceding}~$x$ plus one.

\goodbreak
When replacing each letter~$x$ in the Lyndon word
factorization of~$w$ by~$p(x)$, we obtain a {\it cycle
decomposition} of a permutation~$\sigma$. Furthermore, the
cycles start with their minima and when reading the word from left to right the cycle minima are in {\it decreasing} order.

When this replacement is applied to the Lyndon word
factorization displayed in (3.10), we obtain:
$$
\matrice{w&=&
2&\mid&3&2&1&1&\mid& 3&\mid& 5&\mid&
6&4&2&1&3&2&3&\mid &6&6&3&1&6&6&2&\mid &6\cr
\sigma&=&
\bf16&\mid&12&18&22&21&\mid&\bf 10&\mid&\bf 7&\mid&
4&8&17&20&11&15&9&\mid &2&5&13&19&3&6&14&\mid &\bf1\cr
}
$$

Let $\overline c_{i}:=p^{-1}(i)$ for $i=1,2,\ldots,n$.
As the permutation $\sigma$ is expressed as the product of its
disjoint cycles, we can form the three-row matrix
$$
\matrice{{\rm Id}&=&1&2&\cdots&n\cr
\sigma&=&\sigma(1)&\sigma(2)& \cdots&\sigma(n)\cr
\overline
c&=&\overline c_1&\overline c_2&
\cdots&\overline c_n\cr}
$$
The essential feature is that the word $\overline c$ just
defined is the monotonic {\it nonincreasing} rearrangement of~$w$ and it
has the property that
$$
\sigma(i)>\sigma(i+1)\Rightarrow 
\overline c_i>\overline c_{i+1}.\leqno(5.5)
$$
See {[GeRe93, DeWa93]} for a detailed proof.
The rest of the proof is routine. Let $z=z_1z_2\cdots z_n$ be
the word defined by
$$
z_i:=\#\{j:i\le j\le n-1,\, \sigma(j)>\sigma(j+1)\}.
$$
In other words, $z_i$ is the number of descents of~$\sigma$ within
the right factor $\sigma(i)\sigma(i+1)\cdots\sigma(n)$.
In particular, $z_1=\des\sigma$.
Because of (5.5) the word
$c=c_1c_2\cdots c_n$ defined by $c_i:=\overline c_{i}- z_i$ for
$i=1,2,\ldots,n$ belongs to $\NIW(r-\des\sigma)$ and
$\des\sigma\le r$. 

\goodbreak
Finally, the verification of the three
properties (5.3) is straightforward. Thus, we have constructed
the desired bijection $\psi^{\fix}:w\mapsto (\sigma,c)$, as the reverse construction requires no further development.

With the above example we have:
$$\matrice{
{\rm Id}&=&1&2&3&4&5&6&7&8&9&10&11&12&13
&14&15&16&17&18&19&20&21&22\cr
\sigma&=&\bf1&\underline
5&\underline6&\underline8&\underline{13}&\underline{14}&
\bf7&\underline{17}&4&\bf10&\underline{15}&\underline{18}&
\underline{19}&2&9&\bf16&\underline{20}&\underline{22}
&3&11&12&21\cr
\overline
c&=&6&6&6&6&6&6&5&4&3&3&3&3&3&2&2&2&2&2&1&1&1&1\cr
z&=&4&4&4&4&4&4&3&3&2&2&2&2&2&1&1&1&1&1&0&0&0&0\cr
c&=&2&2&2&2&2&2&2&1&1&1&1&1&1&1&1&1&1&1&1&1&1&1\cr}
$$
The excedances of $\sigma$ have been underlined ($\exc\sigma=11$).
As $\tot z=\maj \sigma$, we have $74=\tot w=\maj\sigma+\tot c
=45+29$.
The fixed points are written in boldface ($\fix\sigma=4$).
\medskip
The bijection $\psi^{\pix}:w\mapsto (\sigma,c)$ of
$W_n(r)$ onto ${\goth S}_n(r,\ides)$ is constructed 
by means of the classical {\it standardisation} of words. 
Read~$w$ from left to right and label 1, 2, \dots\ all the maximal letters. 
If there are~$m$ such letters, restart the reading from left to right and label $m+1$, $m+2$, \dots\ the second greatest letters. Pursue this reading method until reaching the minimal letters. Call $\sigma=\sigma(1)\sigma(2)\cdots\sigma(n)$ the permutation derived by reading those labels from left to right.

The permutation $\sigma$ and the word $w$
have the same hook-factorization {\it type}. This means that if
$ah_1h_2\ldots h_s$ (resp. $bp_1p_2\ldots p_k$) is the hook-factorization
of $\sigma$ (resp. $H$-factorization of $w$), then $k=s$ and 
$\lambda a=\lambda b$.  For each $1\leq i \leq k$ we have
$\lambda h_i=\lambda p_i$  and $\inv(h_i)=\rinv(p_i)$. Hence
$ \wlec w=\lec\sigma$ and $ \wpix w=\pix\sigma$.

Now define the word
$z=z_1z_2\ldots z_n$ as follows. If $\sigma(j)=n$ is the maximal letter, 
then $z_j:=0$; if $\sigma(j)=\sigma(k)-1$ and $j<k$, then $z_j:=z_k$; 
if $\sigma(j)=\sigma(k)-1$ and $j>k$, then $z_j:=z_k+1$.  We can verify that
$\imaj\sigma=\tot z$. With $w=x_{1}x_{2}\cdots x_{n}$ define the word $d=d_{1}d_{2}\cdots d_{n}$  by $d_{i}:=x_{i}-z_{i}$ $(1\le i\le n)$. As $z_{j}=z_{k}+1\Rightarrow x_{j}\ge x_{k}+1$, the letters of~$d$ are all nonnegative. 
The final word $c$ is just defined to be the monotonic
nonincreasing rearrangement of~$d$. Finally, properties (5.4)
are easily verified. 

For defining the reverse of $\psi^{\pix}$ we just have to remember that the following inequality holds: $\sigma(j)<\sigma(k)\Rightarrow d_{j}\ge d_{k}$.
This achieves the proofs of Theorems 1.1 and 1.2.
For example, 
$$
\matrice{
{\rm Id}&=&1&2&3&4&5&6&7&8&9&10&11&12&13&14&15&16&17&18&19&20&21&22\cr
w&=&
6&5&3&2&1&3&2&1&3&6&4&1&2&2&3&1&6&6&3&2&6&6\cr
\sigma&=&
1&7&9&14&19&10&15&20&11&2&8&21&16&17&12&22&3&4&13&18&5&6\cr
z&=&
4&3&2&1&0&2&1&0&2&4&3&0&1&1&2&0&4&4&2&1&4&4\cr
d&=&
2&2&1&1&1&1&1&1&1&2&1&1&1&1&1&1&2&2&1&1&2&2\cr
c&=&
2&2&2&2&2&2&2&1&1&1&1&1&1&1&1&1&1&1&1&1&1&1\cr
}
$$

\medskip
\centerline{\bf 6. A bijection on words and the proof of Theorem 1.3}

\medskip
Consider the two bijections $\phi^{\fix}$ and $\phi^{\pix}$ that have been constructed in sections~3 and~4 and consider the bijection~$F$ defined by the following diagram:

\newbox\boxdiag
\setbox\boxdiag=\vbox{\offinterlineskip
\fleche(0,0)\dir(1,1)\long{10}
\centerput(-7,-2){$D_{n}^*(r)$}
\centerput(17,10){$W_{n}(r)$}
\fleche(17,8)\dir(0,-1)\long{17}
\fleche(0,-2)\dir(1,-1)\long{10}
\centerput(17,-13){$W_{n}(r)$}
\centerput(4,7){$\phi^{\fix}$}
\centerput(3,-11){$\phi^{\pix}$}
\centerput(19,-2){$F$}
}

\vskip.7cm
$$\box\boxdiag\kern2.5cm$$

\vskip1.3cm
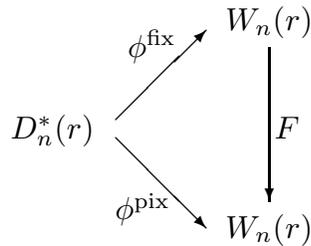
\centerline{Fig. 1}

\bigskip
\noindent
On the other hand, go back to the definition of the bijection $(w,i)\mapsto v$ 
(resp. $(w,i)\mapsto h$) given in Proposition~3.1 (resp. in Proposition~4.1). If $w=x_{1}x_{2}\cdots x_{n}$, then {\it both}~$v$ and~$h$ are {\it rearrangements} of the word $(x_{1}+1)(x_{2}+1)\cdots (x_{i}+1)x_{i+1}\cdots x_{n}$.
Now consider the two bijections
$$\eqalignno{\phi^{\fix}&:(w_{0},(w_{1},i_{1}),\ldots, (w_{m},i_{m}))\mapsto w;\cr
\phi^{\pix}&:(w_{0},(w_{1},i_{1}),\ldots, (w_{m},i_{m}))\mapsto w'.\cr}
$$
It then follows from Proposition~3.2, Proposition~3.3 and (3.9), on the one hand, and from the very definition of~$\phi^{\pix}$, on the other hand, that the words~$w$ and~$w'$ are {\it rearrangements of each other}. Finally, Theorems~2.1 
and 2.3 imply the following result.

\proclaim Theorem 6.1. The transformation $F$ defined by the diagram of~Fig.~$1$ maps each word whose letters are nonnegative integers on another word~$F(w)$ and has the following properties:\hfil\break
\indent(i) $F(w)$ is a rearrangement of $w$ and the restriction of $F$ to each rearrangement class is a bijection of that class onto itself;\hfil\break
\indent(ii) $(\dec,\single)\,w=(\wlec,\wpix)\,F(w)$.

Let {\bf c} be the {\it complement to} $(n+1)$ that maps each permutation
$\sigma=\sigma(1)\sigma(2)\cdots \sigma(n)$ onto
${\bf c}\,\sigma:=(n+1-\sigma(1))(n+1-\sigma(2))\cdots (n+1-\sigma(n))$.
When restricted to the symmetric group ${\goth S}_{n}$ the mapping
$F\circ {\bf c}$ maps ${\goth S}_{n}$ onto itself and has the property
$$
(\des,\single)\,\sigma=(\lec,\pix)\,(F\circ {\bf c})(\sigma).
$$
Note that ``dec" was replaced by ``des", as all the decreases in a permutation are descents. Finally, the so-called first fundamental transformation (see [FoSch70]) $\sigma\mapsto\hat\sigma$ maps ${\goth S}_{n}$ onto itself and is such that
$$\displaylines{
(\exc,\fix)\,\sigma=(\des,\single)\,\hat\sigma.\cr
\noalign{\hbox{Hence}}
(\exc,\fix)\,\sigma=(\lec,\pix)\,(F\circ {\bf c})(\hat\sigma).\cr}
$$
As announced in the introduction we have a stronger result stated in
Theorem~1.3. Its proof is as follows.

\medskip
{\it Proof of Theorem~$1.3$}.\quad
For each  composition $J=j_1j_2\cdots j_m$ (word with positive letters) define the set $L(J)$ and
the monotonic nonincreasing word $c(J)$ by
$$\eqalignno{
L(J)&:=\{j_m, j_m+j_{m-1}, \ldots, j_m+j_{m-1}+\cdots+j_2+j_1\};\cr
c(J)&:=m^{j_m} (m-1)^{j_{m-1}} \cdots 2^{j_2} 1^{j_1}.\cr}
$$
For example, with $J=455116$ we have $L(J)=\{6,7,8,13,18,22\}$ and 
$c(J)=6666665433333222221111$.

Fix a composition $J$ of $n$ (i.e., $\tot J=n$) and let $\Sym^J$ be the set of all permutations $\sigma$ of order $n$
such that $\Iligne\sigma\subset L(J)$.
Using the bijection $\psi^{\pix}$ given in Section~5,
define a bijection $w_1\mapsto \sigma_1$ between   
the set $R_J$ of all rearrangements of $c(J)$ and $\Sym^J$ by
$$
(\sigma_1, *)={\psi^{\pix}} (w_1).
$$
For defining the reverse $\sigma_1\mapsto w_1$ we only have to take
the multiplicity of $w_1\in R_J$ into account.
This is well-defined because $\Iligne\sigma_1\subset L(J)$.
For example, take the same example used in Section~5 for $\psi^{\pix}$:
$$
\matrice{
{\rm Id}&=&1&2&3&4&5&6&7&8&9&10&11&12&13&14&15&16&17&18&19&20&21&22\cr
w_1&=&
6&5&3&2&1&3&2&1&3&6&4&1&2&2&3&1&6&6&3&2&6&6\cr
\sigma_1&=&
1&7&9&14&19&10&15&20&11&2&8&21&16&17&12&22&3&4&13&18&5&6\cr
}
$$
Then
$\Iligne\sigma_1=\{6,8,13,18\}\subset\L(J)$
and the basic properties of this bijection are
$$ \wlec w_1=\lec \sigma_1,\quad \wpix w_1=\pix \sigma_1 .$$

On the other hand the bijection $\psi^{\fix}$ given in Section~5
defines a bijection $w_2\mapsto \sigma_2$ between $R_J$ and $\Sym^J$ by
$$
(\sigma_2^{-1}, *)={\psi^{\fix}} (w_2).
$$
Again, for the reverse $\sigma_2\mapsto w_2$
the multiplicity of $w_2\in R_J$ is to be taken into account.
This is also well-defined, since $\Iligne\sigma_2\subset L(J)$.
With the example used in Section~5 for $\psi^{\fix}$ we have:
$$\matrice{
{\rm Id}&=&1&2&3&4&5&6&7&8&9&10&11&12&13
&14&15&16&17&18&19&20&21&22\cr
w_2&=& 2&3&2&1&1& 3& 5&
6&4&2&1&3&2&3&6&6&3&1&6&6&2&6\cr
\sigma_2^{-1}&=&\bf1& 5&6&8&{13}&{14}& \bf7&{17}&4&\bf10&{15}&{18}&
{19}&2&9&\bf16&{20}&{22} &3&11&12&21\cr
\sigma_2&=&\bf1& 14&19&9&2&3& \bf7&4&15&\bf10&20&21&
5&6&11&\bf16&8&12 &13&17&22&18\cr
}
$$
Also
$\Iligne\sigma_2=\Ligne\sigma_2^{-1}=\{6,8,13,18\}\subset\L(J)$.

\goodbreak
The basic properties of this bijection are
$$ \dec w_2=\iexc \sigma_2,\quad \single w_2=\fix \sigma_2.$$ 
We can use those two bijections and the bijection~$F$ defined in Fig.~1
to form the chain
$$
\sigma\mapsto w_{1}\buildrel F\over \mapsto w_{2}\mapsto \sigma_{2},
$$
and therefore obtain a bijection $\sigma_1\mapsto \sigma_2$ of $\Sym^J$ onto itself having the following properties
$$
\iexc\sigma_2=\lec\sigma_1, \quad
\fix\sigma_2=\pix\sigma_1.
$$
In other words, the pairs $(\iexc, \fix)$ and $(\lec, \pix)$ are equidistributed on
$\{\sigma\in\Sym_n, \Iligne\sigma\subset J\}$ for all compositions $J$ of $n$.
By the inclusion-exclusion principle those pairs are also equidistributed on 
each set $\{\sigma\in\Sym_n, \Iligne\sigma= J\}$. Hence 
the triplets $(\iexc, \fix,\Iligne)$ and $(\lec, \pix, \Iligne)$ are equidistributed on~$\Sym_{n}$.\qed

\bigskip

\centerline{\bf 7. Proof of Theorem 1.4}

\medskip
If ${\bf m}=(m_1,m_2,\ldots,m_n)$ is a sequence of $n$ nonnegative integers, the rearrangement class of the nondecreasing word $1^{m_1}2^{m_2}\ldots r^{m_n}$, that is, the class of all the words than can be derived from
$1^{m_1}2^{m_2}\ldots r^{m_n}$ by permutation of the letters,
is denoted by $R_{\bf m}$. The definitions of ``des," ``maj" and ``inv" 
used so far for permutations are also valid for words. The {\it second fundamental transformation}, as it was called later on (see [Fo68], {[Lo83]}, \S 10.6 or {[Kn73]}, ex. 5.1.1.19) denoted by~$\Phi$, maps each word~$w$ on another
word
$\Phi(w)$ and has the following properties:

(a) $\maj w=\inv \Phi(w)$;

(b) $\Phi(w)$ is a rearrangement of~$w$, and the restriction
of~$\Phi$ to each rearrangement class $R_{\bf m}$ is a bijection
of $R_{\bf m}$ onto itself.

Further properties were further proved by Foata,
Sch\"utzenberger [FoSch78] and Bj\"orner, Wachs [BjW91], in
particular, when the transformation is restricted to act on
rearrangement classes $R_{\bf m}$ such that $m_1=\cdots=m_n=1$,
that is, on symmetric groups~${\goth S}_n$. 

Ligne and inverse ligne of route have been defined in the Introduction. As was proved in
{[FoSch78]}, the transformation $\Phi$ preserves the inverse ligne of route, so that the pairs $(\Iligne,\maj)$ and $(\Iligne,\inv)$ are equidistributed on ${\goth S}_{n}$, a result that we express as
$$\leqalignno{
(\Iligne, \maj)&\simeq(\Iligne, \inv);&(7.1)\cr
\noalign{\hbox{or as}}
(\Ligne, \imaj)&\simeq(\Ligne, \inv).&(7.2)\cr}
$$

\goodbreak
The refinement of (7.2) we now derive (see Proposition 7.1 and Theorem 7.2 below) is based on the properties of a new statistic called $\LAC$.  
\medskip
For each permutation $\sigma=\sigma(1)\sigma(2)\cdots\sigma(n)$ and each integer~$i$ such that $1\le i\le n$ define
$\ell_{i}:=0$ if $\sigma(i)<\sigma(i+1)$ and 
$\ell_{i}:=k$ if $\sigma(i)$ is greater than all the letters
$\sigma(i+1)$, $\sigma(i+2)$, \dots, $\sigma(i+k)$, but
$\sigma(i)<\sigma(i+k+1)$. [By convention, $\sigma(n+1)=+\infty$.]

\medskip
{\it Definition}.\quad
The statistic $\LAC\sigma$ attached to each permutation
$\sigma=\sigma(1)\sigma(2)\cdots\sigma(n)$ is defined to be the word 
$\LAC\sigma=\ell_1\ell_2\ldots \ell_n$.
\medskip
{\it Example}. We have
$$
\matrix{
\hbox{id} &=&1&2&3&4&5&6&7&8&9&10&11&12\cr
\sigma    &=&3&4&8&1&9&2&5&10&12&7&6&11\cr
\LAC\sigma&=&0&0&1&0&2&0&0& 0& 3&1&0& 0\cr
}
$$
In the above table $\ell_5=2$ because 
$\sigma=\ldots {\bf 9}\ 2\ 5\ {\bf 10} \cdots$
and $\ell_9=3$ because
$\sigma=\ldots {\bf 12}\ 7\ 6\ 11$.

\proclaim Proposition 7.1. Let $\sigma=\sigma(1)\sigma(2)\cdots\sigma(n)$ be
a permutation and let $\LAC\sigma=\ell_1\ell_2\ldots \ell_n$.
Then $i\in\Ligne\sigma$ if and only if $\ell_i\geq 1$.

\proclaim Theorem 7.2. We have
$$
(\LAC, \imaj)\simeq(\LAC, \inv).\leqno{(7.3)}
$$

{\it Proof}.\quad
Define $\ILAC\sigma:=\LAC\sigma^{-1}$. Since $\Phi$ maps ``$\maj$" to ``$\inv$," 
property (7.3) will be proved if we show that $\Phi$ preserves ``$\ILAC$", that~is,
$$
\ILAC\Phi(\sigma)=\ILAC\sigma.\leqno(7.4)
$$
A direct description of $\ILAC\sigma$ can be given as follows. Let
$\ILAC\sigma=f_1f_2\ldots f_n$. Then $f_i=j$ if and only if 
within the word $\sigma=\sigma(1)\sigma(2)\cdots\sigma(n)$
the integer $j$ is such that the letters of~$\sigma$ equal to $i+1$, $i+2$, \dots~,  $i+j$ are on the left of the letter equal to~$i$
and either $(i+j+1)$ is on the right of $i$, or $i+j=n$. 

As can be seen in ({[Lo83]}, chap. 10), the second fundamental transformation $\Phi$ is defined by induction: $\Phi(x)=x$ for each letter~$x$ and $\Phi(wx)=\gamma_x(\Phi(w))x$ for each word~$w$ and each letter~$x$, where $\gamma_x$ is a well-defined bijection. See the above reference for an explicit description of $\gamma_x$. Identity (7.4) is then a simple consequence of the following property of~$\gamma_{x}$ (we omit its proof):
{\sl Let $w$ be a word and $x$ a letter. 
If $u$ is a subword of $w$ such that all letters of $u$ are smaller 
(resp. greater) than $x$, then $u$ is also a subword of $\gamma_x(w)$.}\qed

\proclaim Proposition 7.3. Let $\sigma$ and $\tau$ be two permutations
of  order $n$. If $\LAC\sigma=\LAC\tau$, then\hfil\break
\indent(i) $\Ligne\sigma=\Ligne\tau$;\hfil\break
\indent(ii) $(\des,\maj)\sigma=(\des,\maj)\tau$;\hfil\break
\indent(iii) $\pix\sigma=\pix\tau$;\hfil\break
\indent(iv) $\lec\sigma=\lec\tau$.
\medskip

{\it Proof}. 
(i) follows from Proposition 7.1. (ii) follows from (i). By (i) we see that $\sigma$ and $\tau$
have the same hook-factorization {\it type}. That means that if
$ah_1h_2\ldots h_s$ (resp. $bp_1p_2\ldots p_k$) is the hook-factorization
of $\sigma$ (resp. of $\tau$), then $k=s$ and 
$\lambda\,a=\lambda\,b$,  
$\lambda\,h_i=\lambda\,p_i$ for $1\leq i \leq k$. Hence (iii) holds.
For proving (iv) it suffices to prove that  $\inv(h_i)=\inv(p_i)$ for
$1\leq i\leq k$. This is true since $\LAC\sigma=\LAC\tau$
by hypothesis.\qed

\medskip
It follows from Proposition 7.3 that
$$
(\lec, \imaj,\pix)\simeq(\lec,\inv,\pix)\leqno{(7.5)}
$$
and this is all we need to prove Theorem 1.4.

\bigskip\bigskip

\vfill\eject\vglue 1.5cm
\centerline{\bf References}

\bigskip
{\eightpoint

\article BjW91|Anders Bj\"orner, Michelle L. Wachs|Permutation Statistics and Linear Extensions of Posets|J. Combin. Theory, Ser.~A|58|1991|85--114|

\livre An76|George E. Andrews|The Theory of
Partitions|Addison-Wesley, Reading MA, {\oldstyle 1976} ({\it
Encyclopedia of Math.and its  Appl., {\bf 2}})|

\article Ch58|K.T. Chen, R.H. Fox, R.C. Lyndon|Free differential calculus, IV. The quotient group of the lower central series|Ann. of Math.|68|1958|81--95|

\article CHZ97|R. J. Clarke, G.-N. Han, J. Zeng|A combinatorial
interpretation of the Seidel generation of $q$-derangement numbers|%
Annals of Combinatorics|4|1997|313--327|

\divers De84|Jacques D\'esarm\'enien|Une autre interpr\'etation du nombre 
de d\'erangements, {\it 
S\'emi\-naire Lotharin\-gien de Combinatoire}, [B08b], {\oldstyle
1984},  6 pages|

\divers DeWa88|Jacques D\'esarm\'enien, Michelle L.
Wachs|Descentes des d\'erangements et mots circulaires, {\sl
S\'em. Lothar. Combin.}, B19a, {\oldstyle1988}, 9 pp. (Publ.
I.R.M.A. Strasbourg, {\oldstyle1988}, 361/S-19, p.~13-21)|

\article DeWa93|Jacques D\'esarm\'enien, Michelle L.
Wachs|Descent Classes of Permutations with a Given Number of
Fixed Points|J. Combin. Theory, Ser.~A|64|1993|311--328|

\article Fo68|Dominique Foata|On the Netto inversion 
number of a sequence|Proc. Amer. Math. Soc.|19|1968|236--240|

\divers FoHa06a|Dominique Foata, Guo-Niu Han|Fix-Mahonian Calculus, I: two transformations, preprint 16~p., {\oldstyle 2006}|

\divers FoHa06b|Dominique Foata, Guo-Niu Han|Fix-Mahonian Calculus, II: further statistics, preprint 13~p., {\oldstyle 2006}|

\article FoSch78|Dominique Foata,
M.-P. Sch\"utzenberger|Major Index and Inversion
number of Permutations|Math. Nachr.|83|1978|143--159|

\article Ge91|Ira Gessel|A coloring 
problem|Amer. Math. Monthly|98|1991|530--533|

\article GeRe93|Ira Gessel, Christophe Reutenauer|Counting
Permutations with Given Cycle Structure and Descent Set| J.
Combin. Theory Ser. A|64|1993|189--215|


\article KiZe01|Dongsu Kim, Jiang Zeng|A new decomposition of derangements|J. Combin. Theory Ser. A|96|2001|192--198|

\livre Kn73|Donald E. Knuth|The Art of Computer Programming,  {\rm vol. 3}, Sorting and
Searching|Addison-Wesley, Reading
{\oldstyle1973}|

\livre Lo83|M. Lothaire|Combinatorics on Words|Addison-Wesley,
London {\oldstyle 1983} (Encyclopedia of Math. and its Appl., {\bf
17})|

\article Sch65|M.-P. Sch\"utzenberger|On a factorization of free monoids|Proc. Amer. Math. Soc.|16|1965|21-24|

\divers ShWa06|John Shareshian, Michelle L. Wachs|$q$-Eulerian
Polynomials: Excedance Number and Major Index, preprint, 14 p.,
{\oldstyle 2006}|

}

\bigskip\bigskip
\hbox{\vtop{\halign{#\hfil\cr
Dominique Foata \cr
Institut Lothaire\cr
1, rue Murner\cr
F-67000 Strasbourg, France\cr
\noalign{\smallskip}
{\tt foata@math.u-strasbg.fr}\cr}}
\qquad
\vtop{\halign{#\hfil\cr
Guo-Niu Han\cr
I.R.M.A. UMR 7501\cr
Universit\'e Louis Pasteur et CNRS\cr
7, rue Ren\'e-Descartes\cr
F-67084 Strasbourg, France\cr
\noalign{\smallskip}
{\tt guoniu@math.u-strasbg.fr}\cr}}
}

\bye